\documentclass{amsart}

\usepackage{graphicx}
\usepackage{bbm,amsmath,amsfonts,amssymb}
\usepackage[mathscr]{eucal}

\usepackage{natbib}

\newtheorem{theorem}{Theorem}
\newtheorem{lemma}{Lemma}
\newtheorem{corollary}{Corollary}
\newtheorem{proposition}{Proposition}
\newtheorem{definition}{Definition}
\newtheorem{assumption}{Assumption}

\newtheorem{remark}{Remark}

\newcounter{figures}[section]

\RequirePackage[colorlinks,allcolors =black]{hyperref}

\newtheorem{assumption*}{Assumption}
\newtheorem*{assA1}{Condition $\mathbf{A_1}$}
\newtheorem*{assA2}{Condition $\mathbf{A_2}$}
\newtheorem*{assA3}{Condition $\mathbf{A_3}$}

\renewcommand{\kappa}{\varkappa}

\newcommand{\rd}{{\rm d}}
\newcommand{\cD}{{\mathcal D}}
\newcommand{\cE}{{\mathcal E}}
\newcommand{\cF}{{\mathcal F}}

\newcommand{\cI}{{\mathcal I}}

\newcommand{\cM}{{\mathcal M}}

\newcommand{\cR}{{\mathcal R}}

\newcommand{\cX}{{\mathcal X}}
\newcommand{\cY}{{\mathcal Y}}

\newcommand{\Be}{\boldsymbol{e}}

\newcommand{\BL}{\boldsymbol{L}}

\newcommand{\blb}{\boldsymbol{\beta}}

\newcommand{\bga}{\boldsymbol{\gamma}}

\newcommand{\bma}{\boldsymbol{\sigma}}

\newcommand{\bka}{\boldsymbol{\kappa}}
\newcommand{\bbs}{\boldsymbol{s}}
\newcommand{\bbr}{\boldsymbol{r}}
\newcommand{\bA}{\mathbb A}
\newcommand{\bB}{\mathbb B}
\newcommand{\bC}{\mathbb C}

\newcommand{\bE}{\mathbb E}
\newcommand{\bF}{\mathbb F}

\newcommand{\bL}{{\mathbb L}}

\newcommand{\bN}{{\mathbb N}}

\newcommand{\bP}{{\mathbb P}}

\newcommand{\bR}{{\mathbb R}}

\newcommand{\bT}{{\mathbb T}}

\newcommand{\bY}{{\mathbb Y}}

\newcommand{\mA}{\mathfrak{A}}
\newcommand{\mB}{\mathfrak{B}}

\newcommand{\mF}{\mathfrak{F}}
\newcommand{\mS}{\mathfrak{S}}

\newcommand{\mM}{\mathfrak{M}}

\newcommand{\mI}{\mathfrak{I}}

\newcommand{\mh}{\mathfrak{h}}

\newcommand{\mT}{\mathfrak{T}}

\newcommand{\mz}{\mathfrak{z}}

\newcommand{\epr}{\hfill\hbox{\hskip 4pt
                \vrule width 5pt height 6pt depth 1.5pt}\vspace{0.5cm}\par}

\begin{document}

\title[Part I. Lower bounds]
{Adaptive estimation of $\bL_2$-norm of a probability density and related topics I. Lower bounds.}

\author{G. Cleanthous}
\address{Department of Mathematics and Statistics,
National University of Ireland, Maynooth}
\email{galatia.cleanthous@mu.ie}
\author{A. G. Georgiadis}
\address{School of Computer Science and Statistics,
Trinity College of Dublin}
\email{georgiaa@tcd.ie}
\author{O. V. Lepski}
\address{Institut de Math\'ematique de Marseille, CNRS,\\
Aix-Marseille  Universit\'e}
\email{oleg.lepski@univ-amu.fr}

\subjclass[2010]{62G05, 62G20}

\keywords{estimation of functionals,minimax risk,anisotropic Nikolskii's class,minimax adaptive estimation,adaptive rate of convergence}

\begin{abstract}
We deal with the problem of
the adaptive estimation of the $\bL_2$--norm
of a probability density
on $\bR^d$, $d\geq 1$, from independent observations.
The unknown density
is assumed to be uniformly bounded and  to belong to the union of balls  in the isotropic/anisotropic
Nikolskii's spaces. We will show that the optimally  adaptive estimators over the collection of considered functional classes do no exist.
Also, in the framework of an abstract density model  we present several generic lower bounds related to the adaptive estimation of an arbitrary functional of a probability density. These results having independent interest have no analogue in the existing literature.
In the companion paper \cite{G+N+O-Part2} we prove that established lower bounds are tight and provide with explicit construction of adaptive estimators of $\bL_2$--norm of the density.
\end{abstract}

\date{May, 2023}

\maketitle

\section{Introduction}

 Let $(\cX,\mB(\cX),\nu)$ be a measurable space, and let $X$ be a
$\cX$-valued random variable whose law has the density $f$ with respect to the   measure $\nu$. Let $\Phi(f)$ be a given functional and let $\mF$ be the set of all probability densities verifying $|\Phi(f)|<\infty$.

\vskip0.1cm

We observe $X^{(n)}=(X_1,\ldots, X_n)$, where $X_i, i=1,\ldots,n$, are independent copies of $X$. The goal is to estimate the functional $\Phi$
and by an estimator we mean $X^{(n)}$-measurable real function.   The accuracy of an estimator $\widetilde{\Phi}_n$
is measured by the quadratic risk
\[
 \cR_n[\widetilde{\Phi}_n, f]:=\Big(\bE_f \big[\widetilde{\Phi}_n-\Phi(f)\big]^2\Big)^{1/2},
\]
where $\bE_f$ denotes  expectation with respect to the probability measure
$\bP_f$ of  observation $X^{(n)}$.  Thus, $X^{(n)}\in\cX^n$, where $(\cX^n,\mB(\cX^n))$ is the product space equipped with product measure $\nu_n=\nu\otimes\cdots\otimes\nu$ and $\bP_f$, $f\in\mF$, has the density $\mathfrak{p}_f$ w.r.t. $\nu_n$ given by
\begin{equation}
\label{eq:model}
\mathfrak{p}_f(x)=\prod_{i=1}^n f(x_i),\quad x=(x_1,\ldots,x_n)\in\cX^n.
\end{equation}

\subsection{Minimax theory}
With any
estimator $\widetilde{\Phi}_n$ and any subset $\bF$ of
$\mF$ we associate  {\em the maximal risk} of $\widetilde{\Phi}_n$ on $\bF$:
\begin{equation}
\label{eq:def_maximal_risk}
\cR_n\big[\widetilde{\Phi}_n, \bF\big]:=\sup_{f\in\bF}\cR_n[\widetilde{\Phi}_n, f].
\end{equation}
{\em The minimax risk} on $\bF$ is then
\begin{equation}
\label{eq:def_minimax_risk}
 \cR_n[\bF]:=\inf_{\widetilde{\Phi}_n} \cR_n[\widetilde{\Phi}_n, \bF],
\end{equation}
where $\inf$ is taken over all possible estimators. An estimator $\widehat{\Phi}_n$ is called
{\em optimal in order} or {\em rate--optimal} if
\[
 \cR_n[\widehat{\Phi}_n; \bF] \asymp \cR_n[\bF],\;\;\;n\to\infty.
\]
The rate at which $\cR_n[\bF]$ converges to zero as $n$ tends to infinity is referred to as
{\em the minimax rate of convergence} denoted throughout the paper by  $\phi_n (\bF)$.
\par

\paragraph*{Minimax estimation of $\bL_2$-norm} Let $(\cX,\mB(\cX))=(\bR^d,\mB(\bR^d))$, $d\geq 1$, and let $\nu$ is the Lebesgue measure on the Borel $\sigma$-algebra $\mB(\bR^d)$.
Furthermore  we will write $\rd x$ instead of $\nu(\rd x)$.

Set for any $1\leq p\leq\infty$
 $$
 \|f\|_p:=\bigg[\int_{\bR^d}|f(x)|^p\rd x\bigg]^{1/p},\;\;p<\infty,\quad\; \|f\|_\infty=\sup_{x\in\bR^d}|f(x)|.
 $$

 Choose $\Phi(f)=\|f\|_2$
 and let $\bF=\bN_{\vec{r}, d}(\vec{\beta}, \vec{L})\cap\bB_q(Q)$, where $\bN_{\vec{r}, d}(\vec{\beta}, \vec{L})$ is an anisotropic Nikolskii's class (precise definition of the functional class is given below) and
 \[
\bB_q(Q):=\big\{f:\bR^d\to\bR: \|f\|_q\leq Q\big\},\; q>1,\;Q>0.
\]
The problem of minimax estimation of $\bL_p$-norms of a density, $1<p<\infty$, on $\bN_{\vec{r}, d}(\vec{\beta}, \vec{L})\cap\bB_q(Q)$ was studied in \cite{GL20a}, \cite{GL20b}.

Introduce the following notation:
\begin{align*}
&
\frac{1}{\beta}:=\sum_{j=1}^d\frac{1}{\beta_j},\quad \frac{1}{\omega}:=\sum_{j=1}^d\frac{1}{\beta_jr_j},\quad
\BL:=\prod_{j=1}^dL_j^{1/\beta_j},
\\
&
\tau(s):=1-\frac{1}{\omega}+\frac{1}{\beta s},\;\;\; s\in[1,\infty].
\end{align*}
It is worth mentioning that quantities $\tau(\cdot)$ appear in embedding theorems for Nikolskii's classes; for details
see
\cite{Nikolskii}. Asymptotic behavior of the minimax risks on $\bN_{\vec{r}, d}(\vec{\beta}, \vec{L})\cap\bB_q(Q)$ is conveniently expressed in terms
of
these parameters.
Define
\begin{gather}
\label{eq:rate-exponent}
\mz:=\mz\big(\vec{\beta},\vec{r},q\big)=\left\{
\begin{array}{clc}
\frac{1}{2}\wedge\frac{1}{\tau(1)},\quad&\tau(2)\geq 1;
\\*[2mm]
\frac{1-2/q}{2-2/q-\tau(q)},\quad&\tau(2)< 1,\;\tau(q)<0;
\\*[2mm]
\frac{1}{2}\wedge\frac{\tau(2)}{\tau(1)},\quad&\tau(2)<1,\;\tau(q)\geq 0.
\end{array}
\right.
\end{gather}
The  following result corresponding to the case  $p=2$ can  be deduced from Theorem 1 in \cite{GL20a}.
For any $\vec{\beta}\in (0,\infty)^d$, $\vec{L}\in(0,\infty)^d$, $\vec{r}\in[1,\infty]^d$ and  $q\geq 2$
\begin{equation}
\label{eq:th_LB}
 \cR_n\big[\bN_{\vec{r},d}\big(\vec{\beta},\vec{L}\big)\cap \bB_q(Q)\big]\asymp n^{-\mz}, \quad n\to\infty.
\end{equation}
The rate optimal estimator of $\bL_p$-norm for an arbitrary integer $p\geq 2$ was constructed in \cite{GL20b} under additional assumptions: $\vec{r}\in[1,p]^d\cup[p,\infty]^d$
and $q\geq 2p-1$.

\begin{remark}
\label{rem1}
It easy to check that the function $q\mapsto \frac{1-2/q}{2-2/q-\tau(q)}$ is monotonically increasing and therefore the fastest rate corresponds to $q=\infty$.
In this case $\mz$ takes the form
$
\mz=\tfrac{\omega}{\omega+1}.
$
Additionally,
$$
\tfrac{1}{2}\geq\tfrac{\tau(2)}{\tau(1)}\quad\Leftrightarrow\quad \tfrac{\tau(\infty)}{2\tau(1)}\leq 0\quad\Leftrightarrow\quad \tau(\infty)\leq 0.
$$
Hence, if $q=\infty$ the exponent $\frac{\tau(2)}{\tau(1)}$ does not appear and the rate exponent $\mz$ has the most simple and closed form.
\begin{gather}
\label{eq:rate-exponent-infty}
\mz^*\big(\vec{\beta},\vec{r}\big)=:\mz\big(\vec{\beta},\vec{r},\infty\big)=\left\{
\begin{array}{clc}
\frac{1}{2}\wedge\frac{1}{\tau(1)},\quad&\tau(2)\geq 1;
\\*[2mm]
\tfrac{1}{2-\tau(\infty)},\quad&\tau(2)< 1,\;\tau(\infty)<0;
\\*[2mm]
\frac{1}{2},\quad&\tau(2)<1,\;\tau(\infty)\geq 0.
\end{array}
\right.
\end{gather}

\end{remark}
\noindent We finish this section with the definition of the Nikolskii functional class.
\paragraph{Anisotropic Nikolskii class} Let $(\Be_1,\ldots,\Be_d)$ denote the canonical basis of $\bR^d$.
 For function $G:\bR^d\to \bR$ and
real number $u\in \bR$
{\em the first order difference operator}
with step size $u$ in direction of variable
$x_j$ is defined by
$\Delta_{u,j}G (x):=G(x+u\Be_j)-G(x),\;j=1,\ldots,d$.
By induction,
the $k$-th order difference operator with step size $u$ in direction of  $x_j$ is
\[
 \Delta_{u,j}^kG(x)= \Delta_{u,j} \Delta_{u,j}^{k-1} G(x) = \sum_{l=1}^k (-1)^{l+k}\tbinom{k}{l}\Delta_{ul,j}G(x).
\]
\begin{definition}
\label{def:nikolskii}
For given  vectors $\vec{\beta}=(\beta_1,\ldots,\beta_d)\in (0,\infty)^d$, $\vec{r}=(r_1,$ $\ldots,r_d)\in [1,\infty]^d$,
 and $\vec{L}=(L_1,\ldots, L_d)\in (0,\infty)^d$ we
say that the function $G:\bR^d\to \bR$ belongs to  anisotropic
Nikolskii's class $\bN_{\vec{r},d}\big(\vec{\beta},\vec{L}\big)$ if
 $\|G\|_{r_j}\leq L_{j}$ for all $j=1,\ldots,d$ and
 there exist natural numbers  $k_j>\beta_j$ such that
\[
 \big\|\Delta_{u,j}^{k_j} G\big\|_{r_j} \leq L_j |u|^{\beta_j},\;\;\;\;
\forall u\in \bR,\;\;\;\forall j=1,\ldots, d.
\]
\end{definition}

\subsection{Minimax adaptive theory}
Let $\big\{\bF_\vartheta,\vartheta\in\Theta\big\}$ be the collection of subsets of $\mF$, where $\vartheta$ is a nuisance parameter which may have very complicated structure. Without further mentioning we will consider only  scales of functional classes for which a minimax on $\bF_\vartheta$ estimator (usually depending on $\vartheta$) exists for any $\vartheta\in\Theta$. The minimax rate of convergence  on $\bF_\vartheta$ will be denoted by $\phi_n(\bF_\vartheta)$.

\vskip0.1cm

The problem of adaptive estimation can be formulated as follows:
\textit{is it possible to construct a single estimator $\widehat{\Phi}_n$
 which is  simultaneously minimax on each class
 $\bF_\vartheta,\;\vartheta\in\Theta$, i.e. such that}
\[
\limsup_{n\to\infty} \phi^{-1}_n(\bF_\vartheta)\cR_n\big[\widehat{\Phi}_n; \bF_\vartheta\big]<\infty,\quad\forall \vartheta\in\Theta?
\]
We refer to this question as \textit{the  problem of minimax adaptive
estimation over  the scale of classes }
$\{\bF_\vartheta,\;\vartheta\in\Theta \}$.
If such estimator exists it is called  {\it optimally-adaptive}.

\smallskip

It is well-known that optimally-adaptive estimators do not always exist, see \cite{lepski90}, \cite{lepski92a}, \cite{EfLow},  \cite{tsyb}, \cite{kluch}, \cite{rebelles1} among some others.
Formally  nonexistence of optimally-adaptive estimator means that  there exist $\vartheta_1,\vartheta_2\in\Theta$ such that
\begin{equation}
 \label{eq1:nonex}
\liminf_{n\to\infty}\inf_{\widetilde{\Phi}_n}\sup_{\vartheta\in\{\vartheta_1,\vartheta_2\}}\phi^{-1}_n(\bF_{\vartheta}) \cR_n\big[\widetilde{\Phi}_n; \bF_{\vartheta}\big]=\infty.
\end{equation}
Indeed, since a minimax estimator on $\bF_{\vartheta}$ exists for any $\vartheta\in\Theta$ we can assert that
\begin{equation*}
 0<\liminf_{n\to\infty}\inf_{\widetilde{\Phi}_n}\phi^{-1}_n(\bF_{\vartheta}) \cR_n\big[\widetilde{\Phi}_n; \bF_{\vartheta}\big]<\infty,\quad \forall \vartheta\in\Theta.
\end{equation*}
The latter result means that the optimal (from the minimax point of view) family of normalizations $\{\phi_n\big(\bF_{\vartheta}\big), \vartheta\in\Theta\big\}$ is attainable for each value $\vartheta$, while (\ref{eq1:nonex}) shows that this family is unattainable by any estimation procedure simultaneously for a couple of elements from $\Theta$. This, in its turn, implies that optimally-adaptive over the scale $\{\bF_{\vartheta}, \vartheta\in\Theta\big\}$ does not exist.

\vskip0.1cm

However, the question of constructing a single estimator for all values of the nuisance  parameter $\vartheta\in\Theta$ remains relevant. Hence, if (\ref{eq1:nonex}) holds we need to find an attainable family of normalization and to prove its optimality in some manner. The realization of this  program dates back to \cite{lepski90} and \cite{lepski92a} where the notion of
 \textit{adaptive rate of convergence} was introduced.
 Nowadays there exist several definitions of adaptive rate of convergence and corresponding to this notion criteria of optimality, see
\cite{lepski91}, \cite{tsyb}, \cite{kluch}, \cite{rebelles1}.

\smallskip

\noindent The  definition of the adaptive rate, which will be used in the paper, is based  on the comparison of normalization families. It is used usually when $\Theta$ is a subset of some Euclidian space. Furthermore,  $\mB^{\ell}, \ell\in\bN^*$, denotes the set of all subsets of   $\bR^\ell$ containing an euclidian ball, while $\mT^{(\ell-1)}(\bB), \bB\in\mB^{\ell}$, denotes the set of all topological manifolds being subsets of $\bB$ such that those dimension is equal to $\ell-1$.


\begin{remark}
\label{rem2}
The criterion of optimality below, Definition \ref{def2:adaptive-rate}, is used for the classes with different rates of convergence. Note however (and this is a typical situation) that not all components of the nuisance parameter $\vartheta$ may ``effect"  the rate.  For example, looking at the result presented in (\ref{eq:th_LB}) one can remark that only $\vec{\beta}$, $\vec{r}$ and $q$ appear in the rate description. The parameters $\vec{L}$ and $Q$ appear only in the constant bounding the minimax risk.

 With this in mind, we will further assume that the set of nuisance parameters $\Theta$ can be expressed as
\begin{equation}
\Theta=\Theta_{(r)}\times \Theta_{(c)}=\{\vartheta=(\theta,\bar{\theta}):\;\theta\in\Theta_{(r)},\bar{\theta}\in\Theta_{(c)}\},
 \end{equation}
where the minimax rate $\phi_n(\bF_{\vartheta})$ depends just on $\theta\in\Theta_{(r)}$, while the constant bounding the minimax risk depends also on $\bar{\theta}\in\Theta_{(c)}$.
Also, without further mentioning, we will assume that $\Theta_{(r)}\in\mB^\ell$ for some fixed $\ell\in\bN^*$.
It is worth noting that the part of nuisance parameter  $\bar{\theta}$ will be always assumed unknown, in other words the adaptation is considered w.r.t. $\vartheta\in\Theta$ and not w.r.t. $\theta\in\Theta_{(r)}$ only.

\end{remark}


 We say that  $\cM=\{\mu_n(\theta), \theta\in\Theta_{(r)}\}$ is an {\it admissible normalization family} if there exists  an estimator  $\widehat{\Phi}_{n,\cM}$ such that
$$
\limsup_{n\to\infty}\mu^{-1}_n(\theta) \cR_n\big[\widehat{\Phi}_{n,\cM}; \bF_{\vartheta}\big]< \infty,\quad \forall \vartheta\in\Theta.
$$
Denote by $\mM$ the set of all admissible normalization families. Set for any couple $\Psi,\cM\in\mM$
$$
\Theta_-(\Psi,\cM)=\Big\{\theta\in\Theta_{(r)}:\; \limsup_{n\to\infty}\tfrac{\psi_n(\theta)}{\mu_n(\theta)}=\infty\Big\}.
$$
In other words, for any $\theta\in\Theta_-(\Psi,\cM)$ the estimator $\widehat{\Phi}_{n,\cM}$ outperforms the estimator $\widehat{\Phi}_{n,\Psi}$.

\begin{definition}[\cite{kluch}]
\label{def2:adaptive-rate}
  A normalization family $\Psi=\{\psi_n(\theta), \theta\in\Theta_{(r)}\big\}$ is called adaptive rate of convergence  over the collection of functional classes $\{\bF_{\vartheta}, \vartheta\in\Theta\big\}\}$ if:
\begin{eqnarray*}
 \label{eq1:admiss}
  &(i)&\qquad \Psi\in\mM;
\\*[2mm]
  &(ii)&\qquad
\forall\cM\in\mM\;\;\;\exists \bT\in\mT^{\ell-1}\big(\Theta_{(r)}\big):\;\; \Theta_-(\Psi,\cM)\subset \bT;
\\*[2mm]
 \label{eq2:admiss}
&(iii)& \forall  \theta\in\Theta_-(\Psi,\cM),\; \text{there exists}\; \Theta[\theta]\subset\Theta_{(r)}\; \text{such that}\;
\exists \bB\in\mB^{\ell}:\; \Theta[\theta]\supset \bB\;\;\text{and}\;\;
\nonumber
\\*[2mm]
&&\qquad\;\limsup_{n\to\infty}\tfrac{\psi_n(\theta^{\ast})}{\mu_n(\theta^{\ast})}\tfrac{\psi_n(\theta)}{\mu_n(\theta)}=0,\;\;\;\forall \theta^{\ast}\in\Theta[\theta].
\end{eqnarray*}

\noindent The estimator $\boldsymbol{\widehat{\Phi}_{n,\Psi}}$  is called  adaptive and the quantity
$\boldsymbol{\sup_{\theta\in\Theta_{(r)}}[\psi_n(\theta)\big/\phi_n(\bF_\vartheta)]}$
 is called  price to pay for adaptation.

\end{definition}

First we remark that the existence of optimally-adaptive estimator is reduced  to the admissibility of the normalization family $\varphi=\{\phi_n(\bF_\vartheta), \vartheta\in\Theta\big\}$. Indeed, by definition of minimax rate $\Theta_-(\varphi,\cM)=\emptyset$ for any $\cM\in\mM$.
\vskip0.1cm
Next, the definition above has a very simple informal interpretation. As we already mentioned, in the case of nonexistence of optimally-adaptive estimators {\it any estimation procedure} can be outperformed by another one at least for one value of $\vartheta$. Definition \ref{def2:adaptive-rate} stipulates that an adaptive estimator can be outperformed by some estimator but

\begin{itemize}
\item
on a ``negligible" set $\Theta_-(\Psi,\cdot)$ only;
\item
\smallskip

for any value of the parameter at which this procedure is more accurate  than the adaptive estimator, one finds a parameter set of full dimension  on which its loss is greater in order.

\end{itemize}

\begin{remark}[Uniqueness in order of the adaptive rate]
\label{rem3}
As it was noted in \cite{kluch} the adaptive rate of convergence is unique in order.
More precisely if normalization families  $\Psi$ and $\Psi^*$ are both adaptive rates
of convergence then necessarily
\begin{equation}
\label{eq:uniqueness}
0<\liminf_{n\to\infty}\tfrac{\psi^*_n(\theta)}{\psi_n(\theta)}\leq \limsup_{n\to\infty}\tfrac{\psi^*_n(\theta)}{\psi_n(\theta)}<\infty, \quad \forall \theta\in\Theta_{(r)}.
\end{equation}
Indeed, assume that there exists $\theta^*\in \Theta_{(r)}$ such that $\tfrac{\psi^*_n(\theta^*)}{\psi_n(\theta^*)}\to 0$ when $n\to\infty$.
Since $\theta^*\in \Theta_-(\Psi,\Psi^*)$  and $\Psi$ is an adaptive rate there exist $\bB\in\mB^\ell$ such that $\Theta[\theta^*]\supset \bB$.

On the other hand $\tfrac{\psi^*_n(\theta)}{\psi_n(\theta)}\to \infty$ for any $\theta\in \Theta[\theta^*]$ and, therefore, $\Theta[\theta^*]\subset \Theta_-(\Psi^*,\Psi)$. Since $\Psi^*$ is an adaptive rate there exist $\bT\in\mT^{\ell-1}(\Theta_{(r)})$ such that $\Theta[\theta^*]\subset \bT$.
The obtained contradiction proves (\ref{eq:uniqueness}).
\epr

\end{remark}

\subsection{Objectives and organization of the paper} The present paper is the first part of the research project devoted to the adaptive estimation of a general functional of a density. Our goal is two-fold.

$\mathbf{1^0.}\;$ We start with the adaptive estimation of $\bL_2$-norm of a density over the scale of functional classes
$
\cF_\vartheta=\bN_{\vec{r}, d}(\vec{\beta}, \vec{L})\cap\bB_q(Q),\;\vartheta=\big(\vec{\beta},\vec{r},\vec{L},q,Q\big).
$
We prove that the optimally adaptive estimators of $\bL_2$-norm  do not exist. Additionally we establish several lower bounds for adaptive risks over collection of considered classes for different regimes of rate of convergence.  These results form Section \ref{sec:nonex-L_2-norm}. In the companion paper \cite{G+N+O-Part2} we show that established lower bounds are tight and lead to full description of the adaptive rate of convergence.
It is worth mentioning that the price to pay for adaptation increases to infinity proportionally to $\sqrt[4]{\ln(n)}$, $n\to\infty$.

\vskip0.1cm

$\mathbf{2^0.}\;$ We develop a general theory of adaptive estimation mostly related to the construction of lower bounds in the frameworks of an abstract density model introduced in Sections 1.1-1.2. The first one, given in Theorem \ref{th:generic-adaptive} allows to check the nonexistence of optimally-adaptive estimators of an arbitrary functional $\Phi$. In particular, this result  is milestone for the proof of the results from Section \ref{sec:nonex-L_2-norm}. The result presented in Theorem \ref{th:adaptive-rate} provides the conditions under which a given normalization family is an adaptive rate of convergence. The  results mentioned above forms Section \ref{sec:generic_lb} and to best of our knowledge  they do not have the analogues in the existing literature. Also in Section \ref{sec:isotropic-case} the developed approach is applied to the adaptive estimation of $\bL_2$-norm of a density over the scale of isotropic Nikolskii balls of bounded functions.

 \vskip0.1cm

$\mathbf{3^0.}\;$ The proof of main results are given in Section 4 while technical propositions are postponed to Appendix section.

\section{Nonexistence of optimally-adaptive estimators of $\bL_2$-norm of a density}
\label{sec:nonex-L_2-norm}

Set $\Theta_{(r)}=(0,\infty)^d\times[1,\infty]^d\times [2,\infty]$ and let $\theta=\big(\vec{\beta}, \vec{r}, q\big)\in\Theta_{(r)}$.
Set also $\Theta_{(c)}=(0,\infty)^d\times(0,\infty)$ and let $\bar{\theta}=(\vec{L},Q)\in\Theta_{(c)}$. Put at last $\vartheta=(\theta,\bar{\theta})=\big(\vec{\beta}, \vec{r}, q, \vec{L},Q\big)$.

\smallskip

Furthermore, in order to distinguish the estimation of an arbitrary functional $\Phi$ and the estimation of $\bL_2$-norm we will use $\widehat{N}_n, \widetilde{N}_n,\ldots,$
for the estimators of $\bL_2$-norm.
 Our first objective is to prove the following result.
\begin{theorem}
\label{th:nonexistence-1}
The optimally-adaptive estimators for $\bL_2$-norm of a density over collection of functional classes
$$
\cF_\vartheta=\bN_{\vec{r}, d}(\vec{\beta}, \vec{L})\cap\bB_q(Q),\;\;\vartheta=\big(\vec{\beta},\vec{r},q,\vec{L},Q\big)\in\Theta_{(r)}\times\Theta_{(c)}=:\Theta,
$$
do not exist.
\end{theorem}

\noindent The assertion of the theorem is the direct consequence  of much more strong results formulated below. To present them we need to change slightly  the  previous notations, since in the adaptive setting the parameters $\vec{\beta}$ and  $\vec{r}$ are not fixed.

\smallskip

\noindent Set $\varsigma=(\vec{\beta}, \vec{r})$ and redesignate
$\omega$  and $\tau(s)$ introduced in Section 1 by
\begin{align}
\label{eq:new-tau&omega}
 \tfrac{1}{\omega_{\varsigma}}:=\sum_{j=1}^d\tfrac{1}{\beta_jr_j},\quad\;
\tau_{\varsigma}(s):=1-\tfrac{1}{\omega_{\varsigma}}+\tfrac{1}{\beta s},\;\;\; s\in[1,\infty].
\end{align}
Also instead of $\theta=\big(\vec{\beta}, \vec{r}, q\big)\in\Theta_{(r)}$ we will write $\theta=(\varsigma,q)\in\Theta_{(r)}$. Set
$$
\Theta^{\prime}=\big\{\theta\in\Theta_{(r)}:\; \tau_{\varsigma}(2)\geq 1,\;\tau_{\varsigma}(1)>2\big\}.
$$
Recall that $\mz(\theta), \theta\in\Theta_{(r)}$, is given in (\ref{eq:rate-exponent}) and introduce the following notation.
\begin{gather*}
\Theta^{\prime}[\theta]=\big\{\theta^\prime=(\varsigma^\prime,q^\prime)\in\Theta^{\prime}:\; \mz(\theta^\prime)> \mz(\theta)\big\},\quad\theta=(\varsigma,q)\in\Theta^{\prime}.
\end{gather*}

\begin{proposition}
\label{prop:nonexitence-1}
For any $\theta\in\Theta^{\prime}$, $\theta^\prime\in\Theta^{\prime}[\theta]$ and any $\alpha\in\big(\mz(\theta),\mz(\theta^\prime)\big)$
$$
\liminf_{n\to\infty}\inf_{\widetilde{N}_n}\Big\{n^{2\alpha}\cR^2_n\big[\widetilde{N}_n, \cF_{\vartheta^\prime}\big]+\Big(\tfrac{\sqrt{\ln(n)}}{n}\Big)^{-2\mz(\theta)}\cR^2_n\big[\widetilde{N}_n, \cF_{\vartheta}\big]\Big\}>0.
$$
\end{proposition}
\noindent  Set
$ \Theta^{\prime\prime}=\big\{\theta\in\Theta_{(r)}:\; \tau_{\varsigma}(2)< 1,\;\tau_{\varsigma}(q)< 0\big\}$ and define for any $\theta^\prime=(\varsigma^\prime, q^\prime)$, $\varsigma^\prime=(\vec{\gamma},\vec{s})$,
$$
\rho_{\theta,\theta^\prime}=\sum_{l=1}^d\Big[\tfrac{1}{\beta_l}\big(\tfrac{1}{r_l}-\tfrac{1}{q}\big)\Big]
\wedge\Big[\tfrac{1}{\gamma_l}\big(\tfrac{1}{s_l}-\tfrac{1}{q^\prime}\big)\Big].
$$
Introduce the following set of parameters
\begin{gather*}
\Theta^{\prime\prime}[\theta]=\big\{\theta^{\prime}=(\varsigma^{\prime},q^{\prime})\in\Theta^{\prime\prime}:\; \mz(\theta^{\prime})> \mz(\theta),\;
\rho_{\theta,\theta^\prime}\geq 1,\;q^{\prime}\leq q\big\},\quad\theta=(\varsigma,q)\in\Theta^{\prime\prime}.
\end{gather*}

\begin{remark}
\label{rem5}
Note that $\tau_{\varsigma}(q)<0$ and  $\tau_{\varsigma^{\prime}}(q^\prime)<0$ since $\theta,\theta^\prime\in\Theta^{\prime\prime}$.  Hence
in the isotropic case $\gamma_l=\bga$,  $s_l=\bbs$,  $\beta_l=\blb$ and   $r_l=\bbr$, for any $ l=1,\ldots,d,$ the restriction $\rho_{\theta,\theta^\prime}\geq 1$ appeared in the definition of the set $\Theta^{\prime\prime}[\theta]$ is automatically fulfilled.

\end{remark}

\begin{proposition}
\label{prop:nonexitence-2}
For any $\theta\in\Theta^{\prime\prime}$, $\theta^{\prime}\in\Theta^{\prime\prime}[\theta]$ and any $\alpha\in\big(\mz(\theta),\mz(\theta^{\prime})\big)$
$$
\liminf_{n\to\infty}\inf_{\widetilde{N}_n}\Big\{n^{2\alpha}\cR^2_n\big[\widetilde{N}_n, \cF_{\vartheta^{\prime}}\big]+\Big(\tfrac{\sqrt{\ln(n)}}{n}\Big)^{-2\mz(\theta)}\cR^2_n\big[\widetilde{N}_n, \cF_{\vartheta}\big]\Big\}>0.
$$
\end{proposition}
\noindent

The results of Propositions \ref{prop:nonexitence-1} and \ref{prop:nonexitence-2} will be milestone for establishing of adaptive lower bounds in Section \ref{sec:adaptation_lb-bounded_case}. Their proofs (postponed to "Supplemental materials"), in its turn, are the consequence of very general result given in Theorem \ref{th:generic-adaptive} formulated in Section \ref{sec:generic_lb}.
The following statement is the immediate consequence of  Propositions \ref{prop:nonexitence-1} and \ref{prop:nonexitence-2}.

\begin{corollary}
\label{cor:cor1}
The optimally-adaptive estimators for $\bL_2$-norm of a density exist neither over collection of functional classes
$
\{\cF_\vartheta,\vartheta\in\Theta^{\prime}\}
$
nor over collection
$
\{\cF_\vartheta,\vartheta\in\Theta^{\prime\prime}\}.
$
Additionally,
\begin{equation}
 \label{eq1:nonex1}
\liminf_{n\to\infty}\inf_{\widetilde{N}_n}\sup_{\vartheta\in\{\vartheta_1,\vartheta_2\}}\phi^{-1}_n(\cF_{\vartheta}) \cR_n\big[\widetilde{N}_n; \cF_{\vartheta}\big]=\infty,\quad \forall \vartheta_1,\vartheta_2:\; \theta_1,\theta_2\in\Theta^\prime
\end{equation}
\vspace{2mm}
\begin{equation}
\label{eq1:nonex2}
\liminf_{n\to\infty}\inf_{\widetilde{N}_n}\sup_{\vartheta\in\{\vartheta_1,\vartheta_2\}}\phi^{-1}_n(\cF_{\vartheta}) \cR_n\big[\widetilde{N}_n; \cF_{\vartheta}\big]=\infty,\quad \forall \vartheta_1,\vartheta_2:\; \theta_1\in\Theta^{\prime\prime},\theta_2\in\Theta^{\prime\prime}[\theta_1].
\end{equation}

\end{corollary}

We remark that (\ref{eq1:nonex1}) and (\ref{eq1:nonex2}) are much more stronger than (\ref{eq1:nonex}). Moreover Corollary \ref{cor:cor1} implies obviously the assertion of Theorem \ref{th:nonexistence-1}.

\paragraph*{Proof of Corollary \ref{cor:cor1}} The proof is elementary. Indeed, noting that $\phi_n(\cF_\vartheta)=n^{-\mz(\theta)}$ and, therefore,
$$
\lim_{n\to\infty}\phi_n(\cF_{\vartheta^\prime})n^{\alpha}=0, \quad \lim_{n\to\infty}\phi_n(\cF_{\vartheta})\Big(\tfrac{\sqrt{\ln(n)}}{n}\Big)^{-\mz(\theta)}=0,
$$
and choosing $\vartheta_1=\vartheta^\prime$ and $\vartheta_2=\vartheta$ we deduce (\ref{eq1:nonex1}) from Proposition \ref{prop:nonexitence-1}.
Similarly, (\ref{eq1:nonex2}) is deduced from Proposition \ref{prop:nonexitence-2}.
\epr

\paragraph{\bf Open problem 1} Set
$
\Theta^{\prime\prime\prime}=\big\{\theta\in\Theta_{(r)}:\; \tau_{\varsigma}(2)< 1,\;\tau_{\varsigma}(q)\geq 0,\; \tau_{\varsigma}(\infty)< 0 \big\}
$
and note that
$
\Theta^{\prime}\cup\Theta^{\prime\prime}\cup\Theta^{\prime\prime\prime}
$
is the parameter set on which the minimax rate of convergence of $\bL_2$-norm of a density is nonparametric one. The following problem remains unresolved: the existence of optimally-adaptive estimators for the $\bL_2$-norm of a density
$
\{\cF_\vartheta,\vartheta\in\Theta^{\prime\prime\prime}\}.
$

\paragraph{\bf Open problem 2} The challenging problem, also remained unresolved, is the adaptation with respect to the parameter $q$. More precisely, let $\vec{\beta}, \vec{r}$, $\vec{L}$ and $Q$  be given and consider the adaptive estimation over the scale
$$
\cF_\vartheta=\bN_{\vec{r}, d}(\vec{\beta}, \vec{L})\cap\bB_q(Q),\;\;\vartheta=q.
$$
The minimax rate of convergence depends on $q$ only if $\tau_{\varsigma}(2)<1$ and $\tau_{\varsigma}(q)<0$. We note that Proposition \ref{prop:nonexitence-2} is not applicable in this case because the set $\Theta^{\prime\prime}[\theta]$ is empty. Indeed, if $\vec{\beta}, \vec{r}$ are fixed the inequality $\mz(\theta^{\prime})> \mz(\theta)$ is possible only if $q^\prime>q$.

\paragraph{\bf Open problem 3} Does the result formulated in Theorem \ref{th:nonexistence-1} remain valid when the estimation of $\bL_p$-norm, $1<p<\infty$
is considered, at least in the case of integer $p>2$?  It is worth saying that in the white gaussian noise model when $p$ is non even integer the answer on aforementioned question is negative. In other words optimally-adaptive estimator exists, see \cite{han&co1}.

\paragraph{Some additional references} The topic closely related to our study is the minimax and minimax adaptive estimation of quadratic functionals (or, more generally, smooth functionals). The different results and discussions can be found in \cite{bickel-ritov, Nem90, Donoho-Nuss, birge-massart, picard, beatrice, beatrice-1, LNS, gay-trib, beatrice-2, trib, cai-low, klem, waart, gine-nickl}. \cite{EfLow} is seemed to be the first paper where nonexistence of optimally-adaptive estimation of quadratic functionals was proved. It is worth noting that the majority of mentioned above paper dealt with the univariate statistical models defined on a bounded domain of $\bR$. As it was remarked in \cite{GL20b} the minimax estimation of the norm of a probability density on whole $\bR^d$ over anisotropic functional  classes brings a lot of new effects. And the adaptive estimation of $\bL_2$-norm in the density model is not an exception.

\section{Generic adaptive lower bound in the abstract density model}
\label{sec:generic_lb}

In this section we consider the estimation of an arbitrary functional $\Phi$ in the frameworks of the abstract density model (\ref{eq:model}) under the risk given in (\ref{eq:def_minimax_risk}).
 Our first  goal  is to prove an "abstract" analog of Propositions \ref{prop:nonexitence-1} and \ref{prop:nonexitence-2}. Moreover we find rather general conditions under which a normalization family is an adaptive rate of convergence.

\subsection{Case of two sets} Let $\bF$ and $\cF$ be two given subsets of $\mF$ and let $\phi_n(\bF)$ and $\phi_n(\cF)$ be the minimax rates of convergence under quadratic risk (\ref{eq:def_minimax_risk})  on
$\bF$ and $\cF$ respectively.

Assume that there exists $0<c<1$ such that
\begin{equation}
\label{eq:rate_ratio}
\alpha_n(c):=\frac{\phi_n(\cF)}{[\phi_n(\bF)]^c}\to\infty,\quad n\to\infty.
\end{equation}

\paragraph*{Parameterized family of functions} All quantities introduced below may depend (and usually depend) on the number of observation $n$ but we will omit this dependence in our notations.

\vskip0.1cm

Let $f_0\in\mF$  be  given density.
 Let $M\in\bN^*$ be given integer and let $\Lambda_m:\cX\to\bR$, $m=1,\ldots,M$, be a collection of  measurable square-integrable functions on $\cX$.

Let $\bY\subseteq\bR$ be a given interval.  For any $y=(y_1,\ldots,y_M)\in\bY^{M}$ and $x\in\cX$ set
\begin{equation}
\label{eq:family}
f_{y}(x)=\big(1-\rho_y\big)f_0(x)+\sum_{m=1}^My_m\Lambda_m(x),
\end{equation}
where we have put
$$
\rho_y=\sum_{m=1}^M y_m\lambda_m,\quad \lambda_m=\int_{\cX} \Lambda_m(x)\nu(\rd x).
$$
We remark that since $f_0$ is a probability density, one has for all  $y\in\bY^M$
\begin{equation}
\label{eq:int=1}
\int_{\cX}f_{y}(x) \nu(\rd x)=1.
\end{equation}
Note however that $f_{y}$ is not necessarily a probability density, since it may take negative values.

\paragraph*{Conditions} We will apply the construction given in (\ref{eq:family}) under several conditions imposed on its elements.

\begin{assumption}
\label{ass:func_lambda}
There exists a collection $\{\Pi_m\in\mB(\cX), m=1,\ldots, M\}$  of pairwise disjoint subsets of $\cX$ such that
$$
\text{supp}(\Lambda_m)=\Pi_m,\quad \forall m=1,\ldots,M.
$$
\end{assumption}

\begin{assumption}
\label{ass:density_0}
 $f_0\in\bF$ and $f_0(x)>0$ for all $x\in \cup_{m=1}^M\Pi_m.$

\end{assumption}

\vskip0.1cm

\noindent Let $\bP$ be a probability measure on $\bY$ and let us equip $\bY^M$ by product measure $\bP^M=\bP\times\cdots\times\bP$.

\vskip0.1cm

\begin{assumption}
\label{ass:density_xi_in_class}
$
\bP^M\big(\bY_{\cF})\geq \frac{230}{231},
$
where  $
\bY_{\cF}=\big\{y\in\bY^M:\; f_{y}\in\cF\big\}.
$

\end{assumption}

\vskip0.1cm

\noindent For any $a>0$ set $\bY_\Phi(a)=\big\{y\in\bY^M:\; |\Phi(f_{y})-\Phi(f_{0})|> 2a\}.$
\vskip0.1cm

\begin{assumption}
\label{ass:on-Phi}
There exists  $\psi_n\to 0, n\to\infty$, verifying
$$
\lim_{n\to\infty}\psi_n\phi^{-1}_n(\cF)=+\infty,
$$
such that for all $n$ large enough
$
\bP^M\big(\bY_\Phi(\psi_n) \big)\geq \frac{230}{231}.
$

\end{assumption}

\noindent We will tacitly assume that all introduced below quantities involving $\bP$ exist. Set
$$
 \Sigma_M=e_1\sum_{m=1}^M\lambda_m,\qquad \mS_M=\sup_{m=1,\ldots,M}|\lambda_m|,\qquad e_1=\int_{\bR}u\bP(du).
$$
Set also
$$
\bY_+=\Big\{y\in\bY^M:\; \inf_{x\in\cX}f_{y}(x)\geq 0\Big\}.
$$
In view of (\ref{eq:int=1}) we can assert that  $\{f_0, f_{y}, y\in\bY_+\}$ is a family of probability densities.
\begin{assumption}
\label{ass:one-of-two}
Either $\lambda_m=0$ for all $m=1,\ldots,M$  and $\bY_+=\bY^M$
or $\bY^M\subseteq [0,1]^M$, $\Lambda_m:\cX\to\bR_+$ for all $m=1,\ldots,M$, and
\begin{equation}
\label{eq:ass-on-parameters}
\Sigma_M\leq 0.25,\qquad\; 256\mS_M^2Mn\leq 1.
\end{equation}
for all $n$ large enough.

\end{assumption}

\noindent Set for any $m=1,\ldots,M$
$$
S_m=\int_{\cX} \tfrac{\Lambda^2_m(x)}{f_0(x)}\nu(\rd x).
$$
and remark that $S_m$ is well-defined in view of Assumption \ref{ass:density_0}.

\begin{theorem}
\label{th:generic-adaptive}
Let Assumptions \ref{ass:func_lambda}--\ref{ass:one-of-two} be fulfilled. Let $c\in(0,1)$ be defined in (\ref{eq:rate_ratio}) and suppose additionally that
\begin{equation}
\label{eq:main-assump}
\lim_{n\to\infty}\alpha^{-2}_n(c)\prod_{m=1}^M\int_{\bR^2}e^{\mathbf{i} nS_m uv}\bP(\rd u)\bP(\rd v)=0,\quad \mathbf{i}=\{1,2\}.
\end{equation}
Then
$$
\liminf_{n\to\infty}\inf_{\widetilde{\Phi}}\Big\{\phi^{-2c}_n(\bF)\cR^2_n\big[\widetilde{\Phi}, \bF\big]+\psi_n^{-2}\cR^2_n\big[\widetilde{\Phi}, \cF\big]\Big\}\geq \tfrac{107}{144e}.
$$
\end{theorem}
\noindent The proof of the theorem is postponed to Section \ref{sec:th2+prop1&2}. Theorem \ref{th:generic-adaptive} implies obviously that there is no optimally-adaptive estimators over the collection $\{\cF,\bF\}$.

\subsection{Adaptive rate of convergence}
\label{sec:adaptative-rate_general_case}

The goal of this section is to show how to use the assertion of Theorem \ref{th:generic-adaptive} for finding of  adaptive rates of convergence.

Let $\big\{\bF_\vartheta,\vartheta\in\Theta\big\}$ be the collection of subsets of $\mF$. Recall that $\Theta=\Theta_{(r)}\times\Theta_{(c)}$, $\vartheta=(\theta,\bar{\theta})$, $\theta\in\Theta_{(r)}$, $\bar{\theta}\in\Theta_{(c)}$ and the minimax rate of convergence $\phi(\bF_\vartheta)$ is completely  determined by $\theta$. Furthermore, for notation convenience we will write $\varphi_n(\theta)$ instead of $\phi(\bF_\vartheta)$. Recall at last, that $\Theta_{(r)}\in\mB^\ell$ for some $\ell\in\bN^*$.
Introduce the following notations. Set for any $\theta\in\Theta_{(r)}$
\begin{gather*}
\Theta^{+}[\theta]=\Big\{\theta^\prime\in\Theta_{(r)}:\; \limsup_{n\to\infty}\tfrac{\varphi_n(\theta^{\prime})}{\varphi_n(\theta)}=0\Big\};
\\[2mm]
\Theta^{0}[\theta]=\Big\{\theta^\prime\in\Theta_{(r)}:\;0<\liminf_{n\to\infty}\tfrac{\varphi_n(\theta^{\prime})}{\varphi_n(\theta)}\leq
\limsup_{n\to\infty}\tfrac{\varphi_n(\theta^{\prime})}{\varphi_n(\theta)}<\infty\Big\}.
\end{gather*}

\begin{assA1}
For any $\theta\in\Theta_{(r)}$ one can find $\bB\in\mB^\ell$ and $\bT\in\mT^{\ell-1}\big(\Theta_{(r)}\big)$ such that
$\Theta^{+}[\theta]\supset \bB$ and $\Theta^{0}[\theta]\subset \bT$.
\end{assA1}

\begin{assA2}
For any $\theta\in\Theta_{(r)}$ and  $\theta^\prime\in\Theta^{+}[\theta]$ there exists $0<c(\theta,\theta^\prime)<1$ such that
$$
\alpha_n\big(\theta,\theta^\prime\big):=\frac{\varphi_n(\theta)}{[\varphi_n(\theta^\prime)]^{c(\theta,\theta^\prime)}}\to\infty,\quad n\to\infty.
$$
\end{assA2}

\noindent Let $\Psi=\{\psi_n(\theta), \theta\in\Theta_{(r)}\}$ be a normalization family  and set
$$
\mathfrak{p}_n=\sup_{\theta\in\Theta_{(r)}}\tfrac{\psi_n(\theta)}{\varphi_n(\theta)}
$$
\begin{assA3}
For any $\theta\in\Theta_{(r)}$ and  $\theta^\prime\in\Theta^{+}[\theta]$
$$
\limsup_{n\to\infty}\mathfrak{p}^2_n\big[\varphi_n(\theta^\prime)\big]^{1-c(\theta,\theta^\prime)}=0.
$$
\end{assA3}

\begin{theorem}
\label{th:adaptive-rate}
Let Conditions $\mathbf{A_1}-\mathbf{A_3}$ hold. Assume also that Assumptions \ref{ass:func_lambda}--\ref{ass:one-of-two}  are fulfilled with $\cF=\bF_{\vartheta}$, $\bF=\bF_{\vartheta^{\prime}}$ and $\psi_n(\theta)$  for all $\vartheta,\vartheta^{\prime}$ such that $\theta\in\Theta_{(r)}$ and $\theta^\prime\in\Theta^{+}[\theta]$. Moreover, let us assume that the condition (\ref{eq:main-assump}) of Theorem \ref{th:generic-adaptive} is verified with $\alpha_n\big(\theta,\theta^\prime\big)$ for any $\vartheta,\vartheta^{\prime}$ such that $\theta\in\Theta_{(r)}$ and $\theta^\prime\in\Theta^{+}[\theta]$.
If the normalization family $\Psi$ is admissible then it is the adaptive rate of convergence and $\mathfrak{p}_n$ is the price to pay for adaptation.
\end{theorem}

\begin{remark}
\label{rem4}
It is obvious that the application of Theorem \ref{th:adaptive-rate} requires to construct the parameterized family of functions given in (\ref{eq:family})
for each value of $\vartheta\in\Theta$. In particular, all quantities involved in this construction may depend (and usually depend) on $\vartheta$.

\end{remark}

\paragraph*{Proof of Theorem \ref{th:adaptive-rate}} Since Assumptions \ref{ass:func_lambda}--\ref{ass:one-of-two} together with (\ref{eq:main-assump}) are fulfilled with $\cF=\bF_{\vartheta}$ and $\bF=\bF_{\vartheta^{\prime}}$  for all $\vartheta,\vartheta^{\prime}$ such that $\theta\in\Theta_{(r)}$ and $\theta^\prime\in\Theta^{+}[\theta]$, we deduce  from Theorem \ref{th:generic-adaptive} under  Condition $\mathbf{A_2}$ the following bound. For any $\theta\in\Theta_{(r)}$ and $\theta^\prime\in\Theta^{+}[\theta]$ one has
\begin{equation}
\label{eq1:proof_Th_ad-rate}
\liminf_{n\to\infty}\inf_{\widetilde{\Phi}}\Big\{[\varphi_n(\theta^\prime)]^{-2c(\theta,\theta^\prime)}\cR^2_n\big[\widetilde{\Phi}, \bF_{\vartheta^\prime}\big]+\psi_n^{-2}(\theta)\cR^2_n\big[\widetilde{\Phi}, \bF_\vartheta\big]\Big\}\geq \tfrac{107}{144e}.
\end{equation}
Note also that for any $\theta\in\Theta_{(r)}$ and $\theta^\prime\in\Theta^{+}[\theta]$
\begin{equation}
\label{eq2:proof_Th_ad-rate}
\limsup_{n\to\infty}\tfrac{\psi_n(\theta^{\prime})}{\big[\varphi_n(\theta^\prime)\big]^{c(\theta,\theta^\prime)}}=0
\end{equation}
in view of Condition $\mathbf{A_3}$. Here we also took into account that $\mathfrak{p}_n>1$ since otherwise $\Psi$ is the family of rates of convergence
and the assertion of the theorem follows.

\vskip0.1cm

\noindent Let $\cM=\{\mu_n(\theta), \theta\in\Theta_{(r)}\}$ be an arbitrary  admissible normalization family and recall that
$$
\Theta_-(\Psi,\cM)=\Big\{\theta\in\Theta_{(r)}:\; \limsup_{n\to\infty}\tfrac{\psi_n(\theta)}{\mu_n(\theta)}=\infty\Big\}.
$$
\quad $\mathbf{1^0.}\;$  Let us prove that there exists $\bT\in\mT^{\ell-1}\big(\Theta_{(r)}\big)$ such that
\begin{equation}
\label{eq3:proof_Th_ad-rate}
\Theta_-(\Psi,\cM)\subset \bT.
\end{equation}
If $\Theta_-(\Psi,\cM)=\emptyset$ then (\ref{eq3:proof_Th_ad-rate}) holds. Assume now $\Theta_-(\Psi,\cM)\neq\emptyset$ and let $\theta_*\in \Theta_-(\Psi,\cM)$.

First, let us suppose that there exists $\theta^*\in\Theta_-(\Psi,\cM)$ such that either $\theta_*\in\Theta^{+}[\theta^*]$ or $\theta^*\in\Theta^{+}[\theta_*]$.
Furthermore,  without loss of generality we will assume that $\theta^*\in\Theta^{+}[\theta_*]$.

\vskip0.1cm

Since $\theta_*,\theta^*\in\Theta_-(\Psi,\cM)$ we have
\begin{equation}
\label{eq4:proof_Th_ad-rate}
\limsup_{n\to\infty}\tfrac{\mu_n(\theta_{*})}{\psi_n(\theta_{*})}=\limsup_{n\to\infty}\tfrac{\mu_n(\theta^{*})}{\psi_n(\theta^{*})}=0.
\end{equation}
Moreover, since $\theta^*\in\Theta^{+}[\theta_*]$, we obtain from (\ref{eq4:proof_Th_ad-rate}) and  (\ref{eq2:proof_Th_ad-rate}) applied with $\theta=\theta_*$ and $\theta^\prime=\theta^*$
\begin{equation}
\label{eq5:proof_Th_ad-rate}
\limsup_{n\to\infty}\tfrac{\mu_n(\theta^{*})}{\big[\varphi_n(\theta^*)\big]^{c(\theta_*,\theta^*)}}=\limsup_{n\to\infty}
\tfrac{\psi_n(\theta^{*})}{\big[\varphi_n(\theta^*)\big]^{c(\theta_*,\theta^*)}}\;\tfrac{\mu_n(\theta^{*})}{\psi_n(\theta^{*})}=0.
\end{equation}
Set at last,
$$
\kappa_n(\theta_*,\theta^*)=\bigg\{\tfrac{\big[\varphi_n(\theta^*)\big]^{c(\theta_*,\theta^*)}}{\mu_n(\theta^{*})}\bigg\}\bigwedge \tfrac{\psi_n(\theta_{*})}{\mu_n(\theta_{*})}
$$
and note that in view of (\ref{eq4:proof_Th_ad-rate}) and  (\ref{eq5:proof_Th_ad-rate})
$
\liminf_{n\to\infty}\kappa_n(\theta_*,\theta^*)=\infty.
$
Since $\cM$ is admissible there exists an estimator $\widehat{\Phi}_{n,\cM}$ such that
$$
\limsup_{n\to\infty}\mu^{-1}_n(\theta) \cR_n\big[\widehat{\Phi}_{n,\cM}; \bF_{\vartheta}\big]=:C(\vartheta)< \infty,\quad \forall \vartheta\in\Theta.
$$
Hence, we obtain using (\ref{eq1:proof_Th_ad-rate})
\begin{align*}
C^2&(\vartheta_*)+C^2(\vartheta^*)\geq\limsup_{n\to\infty}\Big\{\mu^{-2}_n(\theta^*) \cR_n^2\big[\widehat{\Phi}_{n,\cM}; \bF_{\vartheta^*}\big]+\mu^{-2}_n(\theta_*) \cR_n^2\big[\widehat{\Phi}_{n,\cM}; \bF_{\vartheta_*}\big]\Big\}
\\*[2mm]
&\geq\limsup_{n\to\infty}\kappa_n^2(\theta_*,\theta^*)\Big\{[\varphi_n(\theta^*)]^{-2c(\theta_*,\theta^*)} \cR_n^2\big[\widehat{\Phi}_{n,\cM}; \bF_{\vartheta^*}\big]+\psi_n^{-2}(\theta_*)\cR_n^2\big[\widehat{\Phi}_{n,\cM}; \bF_{\vartheta_*}\big]\Big\}
\\*[2mm]
&\geq\limsup_{n\to\infty}\kappa_n^2(\theta_*,\theta^*)\inf_{\widetilde{\Phi}}\Big\{[\varphi_n(\theta^*)]^{-2c(\theta_*,\theta^*)}\cR^2_n\big[\widetilde{\Phi}, \bF_{\vartheta^*}\big]+\psi_n^{-2}(\theta_*)\cR^2_n\big[\widetilde{\Phi}, \bF_{\vartheta_*}\big]\Big\}=\infty.
\end{align*}
The obtained contradiction shows that either $\text{card}(\Theta_-(\Psi,\cM))= 1$ and  therefore (\ref{eq3:proof_Th_ad-rate}) holds, or
$$
0<\liminf_{n\to\infty}\tfrac{\varphi_n(\theta)}{\varphi_n(\theta_*)}\leq
\limsup_{n\to\infty}\tfrac{\varphi_n(\theta)}{\varphi_n(\theta_*)}<\infty,\quad\forall \theta\in \Theta_-(\Psi,\cM).
$$
But then  $\Theta_-(\Psi,\cM)\subseteq \Theta^{0}[\theta_*]$ and (\ref{eq3:proof_Th_ad-rate}) follows from Condition $\mathbf{A_1}$. Thus, (\ref{eq3:proof_Th_ad-rate}) is established that means that the assumption (\ref{eq1:admiss}) of Definition \ref{def2:adaptive-rate} is checked.

\vskip0.1cm

$\mathbf{2^0.}\;$ Let $\theta$ be an arbitrary element of $\Theta_-(\Psi,\cM)$. Since $\cM$ is an admissible family one has
$$
\limsup_{n\to\infty}\psi^{-2}_n(\theta) \cR^2_n\big[\widehat{\Phi}_{n,\cM}; \bF_{\vartheta}\big]=0
$$
and we obtain from (\ref{eq1:proof_Th_ad-rate}) that for any $\theta^{\prime}\in\Theta^{+}[\theta]$
$$
\liminf_{n\to\infty}[\varphi_n(\theta^\prime)]^{-2c(\theta,\theta^\prime)} \cR^2_n\big[\widehat{\Phi}_{n,\cM}; \bF_{\vartheta^\prime}\big]\geq \tfrac{107}{144e}.
$$
It implies in its turn that
\begin{equation}
\label{eq6:proof_Th_ad-rate}
\limsup_{n\to\infty}\mu^{-1}_n(\theta^\prime)[\varphi_n(\theta^\prime)]^{c(\theta,\theta^\prime)}\leq C(\vartheta^\prime)\sqrt{\tfrac{144e}{107}},\quad \forall \theta^{\prime}\in\Theta^{+}[\theta].
\end{equation}
We have for any $\theta^{\prime}\in\Theta^{+}[\theta]$
\begin{eqnarray*}
\tfrac{\psi_n(\theta^{\prime})}{\mu_n(\theta^{\prime})}\tfrac{\psi_n(\theta)}{\mu_n(\theta)}&=&
\tfrac{\psi_n(\theta^{\prime})}{[\varphi_n(\theta^\prime)]^{c(\theta,\theta^\prime)}}
\tfrac{[\varphi_n(\theta^\prime)]^{c(\theta,\theta^\prime)}}{\mu_n(\theta^{\prime})}
\tfrac{\psi_n(\theta)}{\phi_n(\theta)}\tfrac{\phi_n(\theta)}{\mu_n(\theta)}\leq
\tfrac{\mathfrak{p}_n\psi_n(\theta^{\prime})}{[\varphi_n(\theta^\prime)]^{c(\theta,\theta^\prime)}}
\tfrac{[\varphi_n(\theta^\prime)]^{c(\theta,\theta^\prime)}}{\mu_n(\theta^{\prime})}
\tfrac{\phi_n(\theta)}{\mu_n(\theta)}
\\*[2mm]
&\leq& \mathfrak{p}^2_n\big[\varphi_n(\theta^\prime)\big]^{1-c(\theta,\theta^\prime)}
\tfrac{[\varphi_n(\theta^\prime)]^{c(\theta,\theta^\prime)}}{\mu_n(\theta^{\prime})}
\tfrac{\phi_n(\theta)}{\mu_n(\theta)},
\end{eqnarray*}
where $\phi_n(\theta)$ the minimax rate of convergence on $\bF_{\vartheta}$.

Taking into account that by definition of minimax rate of convergence for any admissible $\cM$ and any $\theta\in\Theta_{(r)}$
$$
\limsup_{n\to\infty}\tfrac{\phi_n(\theta)}{\mu_n(\theta)}<\infty,
$$
we deduce from Condition $\mathbf{A_3}$ and (\ref{eq6:proof_Th_ad-rate}) that for any $\theta\in\Theta_-(\Psi,\cM)$
$$
\limsup_{n\to\infty}\tfrac{\psi_n(\theta^{\prime})}{\mu_n(\theta^{\prime})}\tfrac{\psi_n(\theta)}{\mu_n(\theta)}=0,\quad \forall \theta^{\prime}\in\Theta^{+}[\theta].
$$

It remains to note that in view of  Condition $\mathbf{A_1}$ for any $\theta\in\Theta_{(r)}$ there exists $\bB\in\mB^\ell$ such that $\Theta^{+}[\theta]\supset\bB$. Thus, the assumption (\ref{eq2:admiss}) of Definition \ref{def2:adaptive-rate} is checked as well. This completes the proof of the theorem.
\epr

\subsection{Some generalization} As we have seen the basic step in establishing of adaptive rates is the lower bound for the sum of normalized risks  given in Theorem \ref{th:generic-adaptive}. This, in its turn, requires, see Remark \ref{rem4}, to find the family $\{f_y\}_y\subset\bF_\vartheta$ while the density $f_0$ should belong to $\bF_{\vartheta^\prime}$.  When the estimation problem in the density model\footnote{This problem does not usually appear in the regression context because one can choose $f_0\equiv 0$ and the latter function belongs to all functional classes used in nonparametric statistics. Unfortunately it is not a probability density.} possesses different regimes of minimax rate of convergence this may bring difficulties, if $\vartheta^\prime$ and  $\vartheta$ correspond to different regimes. This is exactly the case when the estimation of $\bL_2$-norm is considered. Looking at (\ref{eq:rate-exponent}) or (\ref{eq:rate-exponent-infty}) we constant there are three (respectively two) different nonparametric regimes. We do not know how to prove the analog of Theorem \ref{th:generic-adaptive} simultaneously for all regimes. Fortunately, the simple result presented below shows that it is not needed. It turns out that it suffices to prove Theorem \ref{th:generic-adaptive} separately for each regime\footnote{See, for instance Proposition \ref{prop:nonexitence-1} and \ref{prop:nonexitence-2}.}.

\smallskip

Let $\Theta_{(r)}=\cup_{k=1}^K \Theta_{(r),k}$, $K\in\bN^*$, where the parameter sets $\Theta_{(r),k}\in\mB^{\ell}, k=1,\ldots, K$ are pairwise disjoint, that is
$\Theta_{(r),k}\cap\Theta_{(r),l}=\emptyset$ for any $l\neq k$.

\vskip0.1cm

Let $\Psi_k=\{\psi_{n,k}(\theta), \theta\in\Theta_{(r),k}\}$, $k=1,\ldots,K$ be the adaptive rate of convergence over the scale of functional classes $\{\bF_{\vartheta},\vartheta\in\Theta_{(r),k}\times\Theta_{(c)}\}$. Define
$$
\psi_n(\theta)=\sum_{k=1}^K \psi_{n,k}(\theta)\mathrm{1}_{\Theta_{(r),k}}(\theta),\quad \theta\in \Theta_{(r)}.
$$
and let $\Psi=\{\psi_{n}(\theta), \theta\in\Theta_{(r)}\}$.

\begin{theorem}
\label{th:adaptive-rate_K}
If $\Psi$ is an admissible normalization family it is the adaptive rate of convergence over the scale of functional classes $\{\bF_{\vartheta},\vartheta\in\Theta_{(r)}\times\Theta_{(c)}\}$.
\end{theorem}

\paragraph*{Proof of Theorem \ref{th:adaptive-rate_K}.} The proof is elementary. Fix some admissible normalization family $\cM=\{\mu_n(\theta), \theta\in\Theta_{(r)}\}$ and set for any $k=1,\ldots, K$
$$ \Theta_{-,k}(\Psi_k,\cM)=\Big\{\theta\in\Theta_{(r),k}:\; \limsup_{n\to\infty}\tfrac{\psi_{n,k}(\theta)}{\mu_n(\theta)}=\infty\Big\}
=\Big\{\theta\in\Theta_{(r),k}:\; \limsup_{n\to\infty}\tfrac{\psi_{n}(\theta)}{\mu_n(\theta)}=\infty\Big\}.
$$
The last equality follows from the definition of $\psi_{n}(\cdot)$. It is obvious that
$$
\Theta_{-}(\Psi,\cM)=\cup_{k=1}^K\Theta_{-,k}(\Psi_k,\cM).
$$
Since $\Psi_k$, $k=1,\ldots, K$, are adaptive rates for any $k$ there exists $\bT_k\in\mT^{l-1}(\Theta_{(r),k})$ such that $\Theta_{-,k}(\Psi_k,\cM)\subset\bT_k$. Moreover $T_k$  $k=1,\ldots, K$, are pairwise disjoint as subsets of pairwise disjoint sets. At last, by its definition the dimension of each $\bT_k$ is $\ell-1$. We have
$$
\Theta_{-}(\Psi,\cM)\subset \cup_{k=1}^K \bT_k=:\bT.
$$
It remains to note that $\bT\in\mT^{l-1}(\Theta_{(r)})$ as the disjoint  union of finite family of topological manifolds of the same dimension. Thus, the condition (\ref{eq1:admiss}) of Definition \ref{def2:adaptive-rate} is checked.

\vskip0.1cm

Fix $\theta\in\Theta_{-}(\Psi,\cM)$ and let $\mathbf{k}\in{1,\ldots, K}$ be such that $\theta\in\Theta_{-,\mathbf{k}}(\Psi_{\mathbf{k}},\cM)$. Since $\Psi_{\mathbf{k}}$
is an adaptive rate there exists $\bB\in\mB^{\ell}$ such that $\Theta[\theta]\supset \bB$ and
$$
0=\limsup_{n\to\infty}\tfrac{\psi_{n,\mathbf{k}}(\theta^{\ast})}{\mu_n(\theta^{\ast})}\tfrac{\psi_{n,\mathbf{k}}(\theta)}{\mu_n(\theta)}=
\limsup_{n\to\infty}\tfrac{\psi_{n}(\theta^{\ast})}{\mu_n(\theta^{\ast})}\tfrac{\psi_{n}(\theta)}{\mu_n(\theta)},\;\;\forall \theta^{\ast}\in\Theta[\theta].
$$
To get the last equality we have used that $\Theta[\theta]\subset\Theta_{(r),\mathbf{k}}$ and the definition of $\psi_n(\cdot)$. Thus, the condition (\ref{eq2:admiss}) of Definition \ref{def2:adaptive-rate} is checked and theorem is proved.
\epr


\subsection{Application to the adaptive estimation of $\bL_2$-norm over collection of isotropic Nikolskii classes of bounded functions}
\label{sec:isotropic-case}

In this section we study the adaptive estimation over the scale of \textit{isotropic} Nikolskii classes of bounded functions. From now on
$$
\beta_{l}=\blb,\; r_l=\bbr,\; L_l=L,\quad\forall l=1,\ldots, d.
$$
Since $\vec{\beta}=(\blb,\ldots,\blb)$, $\vec{r}=(\bbr,\ldots,\bbr)$ the isotropic Nikolskii class will be denoted by $\bN_{\bbr,d}\big(\blb,L\big)$. Moreover we obviously have for any $s\in[1,\infty]$
$$
\tau_{\varsigma}(s)=1-\frac{d}{\blb \bbr}+\frac{d}{\blb s}, \quad \frac{1}{\omega_\varsigma}=\frac{d}{\blb \bbr}, \quad \frac{1}{\beta}=\frac{d}{\blb}.
$$
The adaptation will be studied over the collection of functional classes
$$
\bN_\vartheta:=\bN_{\bbr,d}\big(\blb,L\big)\cap \bB_\infty(Q), \quad \vartheta=(\blb,\bbr,L,Q).
$$
In this section we adopt the following notations.
For any $b\in (0,\infty]$ set $\Theta^b_{(r)}=(0,b]\times[1,\infty]$, $\Theta_{(c)}=(0,\infty)^2$ and let
$\vartheta=(\theta,\bar{\theta})$, where   $\theta=(\blb,\bbr)\in\Theta^b_{(r)}$ and $\bar{\theta}=(L,Q)\in\Theta_{(c)}$.
Thus, $\vartheta\in\Theta=:\Theta^b_{(r)}\times\Theta_{(c)}$.


Note finally that the exponent figurant  in the minimax rate of convergence  given in (\ref{eq:rate-exponent-infty}), becomes in the "isotropic case"
\begin{gather}
\label{eq:rate-exponent-infty-isotropic}
\mz^*(\theta)=\left\{
\begin{array}{cccc}
\frac{\blb\bbr}{\blb\bbr+d(\bbr-1)},\quad& \bbr\geq 2,\; \blb\bbr< d(\bbr-1);
\\*[2mm]
\frac{\blb\bbr}{\blb\bbr+d},\quad&\bbr< 2,\; \blb\bbr< d;
\\*[2mm]
\frac{1}{2},\quad&\bbr\geq 2,\; \blb\bbr\geq d(\bbr-1);
\\*[2mm]
\frac{1}{2},\quad&\bbr< 2,\; \blb\bbr\geq d.
\end{array}
\right.
\end{gather}

\subsection{Adaptive rate of convergence}
\label{sec:adaptation_lb-bounded_case}

Introduce the following sets of parameters.
$$
\Theta^b_{(r),1}=\big\{\theta\in\Theta^b_{(r)}:\;  \bbr\geq 2,\; \blb\bbr< d(\bbr-1)\big\},\quad
\Theta^b_{(r),2}=\big\{\theta\in\Theta^b_{(r)}:\; \bbr< 2,\; \blb\bbr< d\big\}.
$$
We note that  $\Theta^b_{(r),1}\cup\Theta^b_{(r),2}$ is the set of parameters corresponding to nonparametric regime of minimax rate of convergence, that
is $\mz^*(\theta)<\frac{1}{2}$ for any $\theta\in\Theta^b_{(r),1}\cup\Theta^b_{(r),2}$.

\noindent Consider the following normalization family.
$$
\Psi^b=\bigg\{\psi^*_n(\theta):=\Big(\tfrac{\sqrt{\ln(n)}}{n}\Big)^{\mz^*(\theta)},\;\;\theta\in\Theta^b_{(r),1}\cup\Theta^b_{(r),2}\bigg\}.
$$

\begin{theorem}
\label{th:adaptive-lb-bounded}
For any $b\in (0,\infty]$, such that the normalization family $\Psi^b$ is admissible, it is the adaptive rate of convergence over scale of functional classes
$\big\{\bN_\vartheta, \vartheta\in \big[\Theta^b_{(r),1}\cup\Theta^b_{(r),2}\big]\times \Theta_{(c)}\big\}$.
\end{theorem}
\noindent The proof of the theorem is based on application of Propositions \ref{prop:nonexitence-1}, \ref{prop:nonexitence-2} and Theorem \ref{th:adaptive-rate}, \ref{th:adaptive-rate_K}.

\section{Proof Theorems \ref{th:generic-adaptive}, \ref{th:adaptive-lb-bounded} and Propositions \ref{prop:nonexitence-1},\ref{prop:nonexitence-2}}
\label{sec:th2+prop1&2}

\subsection{Proof of Theorem \ref{th:generic-adaptive}} We break the proof on several steps.

\subsubsection{Auxiliary lemmas}

Set  for any $x\in\cX$ and $y\in\bY^M$
$$
f^*_y(x)=\big[1-\Sigma_M\big]e^{-\frac{\rho_y-\Sigma_M}{1-\Sigma_M}}f_0(x)+\sum_{m=1}^My_m\Lambda_m(x).
$$
and let $\bY_\rho=\{y\in\bY^M:\; |\rho_y-\Sigma_M|\leq 12\cD_M\}$, where $\cD_M^2=\sum_{m=1}^M\lambda_m^2$.

\begin{lemma}
\label{lem:bound_for_f_y}
Let Assumptions \ref{ass:func_lambda}, \ref{ass:density_0} and  \ref{ass:one-of-two} be verified. Then for all $n$ large enough
 $$
 f^*_y(x)\geq f_{y}(x)\geq e^{-1/n}f^*_y(x),\quad \forall x\in\cX,\quad y\in\bY_\rho.
 $$

\end{lemma}

\noindent We note that if $\lambda_m=0$ for all $m=1,\ldots,M$ then $f^*_y\equiv f_{y}$ for all $y\in\bY^M$ and the assertion of the lemma is trivial. Note also that in this case
$\bY_\rho=\bY^M$.
\begin{remark}
\label{rem:after-lem:bound_for_f_y}
Since $f^*_y\geq 0$ for any $y\in\bY^M$ in view of Assumption \ref{ass:one-of-two}, we derive from Lemma \ref{lem:bound_for_f_y} and (\ref{eq:int=1}) that $f_y$ is a probability density for any $y\in\bY_\rho$.
\end{remark}

\noindent We will also need the simple probabilistic result formulated in the next lemma.
Let $\zeta_k, k=1,\ldots K$, be i.i.d random variables defined on some probability space $(\Omega,\mA,P)$ and let $E$ denote the expectation corresponding to $P$.
Assume that $P(|\zeta_1|\leq T)=1$ for some $T>0$.
Let $a_k, k=1,\ldots K$ be given sequence of real numbers.
 Set
$$
\chi_K=\sum_{k=1}^K a_k\big(\zeta_k-E\{\zeta_k\}\big),\quad A_K=\sup_{k=1,\ldots,K} |a_k|,\quad D_K=5T^2A^2_KK.
$$
Define finally for any $b>0$ and $J\in\bN^*$
$$
W_{J,K}(b)=E\Big\{\big[1+b\chi_K^2\big]^J\Big\}.
$$

\begin{lemma}
\label{lem:probab}
For any $J,K\in\bN^*$ and $b>0$ verifying
$
2bJD_K\leq 1,
$
one has $W_J(b)\leq 2$.

\end{lemma}

\noindent The proofs of Lemmas \ref{lem:bound_for_f_y} and \ref{lem:probab} are postponed to Appendix section.

\subsubsection{Proof of the theorem.} Introduce $\bY_\rho^*$ as follows. If $\lambda_m=0$ for all $m=1,\ldots, M$ set $\bY^*_\rho=\bY^M$. If not set $\bY^*_\rho=\bY_\rho$, where $\bY_\rho$ is defined in Lemma \ref{lem:bound_for_f_y}.
Also $\bE^M$ will denote the integration w.r.t. $\bP^M$.
Define
$$
\cY=\bY^*_\rho\cap\bY_{\cF}\cap \bY_\Phi(\psi_n).
$$
\quad $\mathbf{1^0a}.\;$ We have for any estimator\footnote{To simplify the notations, we omit the dependence of the estimator  on $n$.} $\widetilde{\Phi}$ applying the Tchebychev inequality
\begin{eqnarray}
\label{eq:ol12}
R_n\big(\widetilde{\Phi}\big)&=:&\phi^{-2c}_n(\bF)\cR^2_n\big[\widetilde{\Phi}, \bF\big]+\psi_n^{-2}\cR^2_n\big[\widetilde{\Phi},\bF\big]
\nonumber\\*[2mm]
&\geq& \alpha^2_n(c)\sup_{f\in\bF}\bP_f\big\{|\widetilde{\Phi}-\Phi(f)|\geq \psi_n\big\}+\sup_{f\in\cF}\bP_f\big\{|\widetilde{\Phi}-\Phi(f)|\geq \psi_n\big\}
\nonumber\\*[2mm]
&\geq&\alpha^2_n(c)\bP_{f_0}\big\{|\widetilde{\Phi}-\Phi(f_0)|\geq \psi_n\big\}+\bE^M\Big[\mathrm{1}_{\cY}(y)\bP_{f_{y}}\big\{|\widetilde{\Phi}-\Phi\big(f_{y}\big)|\geq \psi_n\big\}\Big]
\nonumber\\*[2mm]
&=&\alpha^2_n(c)\bE_{f_0}\Big[\mathrm{1}\big\{|\widetilde{\Phi}-\Phi(f_0)|\geq \psi_n\big\}\Big]+\bE^M\bE_{f_{y}}\Big[\mathrm{1}_{\cY}(y)\mathrm{1}\big\{|\widetilde{\Phi}-\Phi\big(f_{y}\big)|\geq \psi_n\big\}\Big]
\nonumber\\*[2mm]
&=:& I_1+I_2.
\end{eqnarray}
Here we have also used Assumption \ref{ass:density_0} guaranteeing $f_0\in\bF$, Remark \ref{rem:after-lem:bound_for_f_y} and the definition of $\bY_{\cF}$.
\vskip0.1cm
$\mathbf{1^0b}.\;$ Some remarks are in order. First, note that for any $y\in\bY_\Phi(\psi_n)$ one has
$$
\mathrm{1}\big\{|\widetilde{\Phi}-\Phi\big(f_{y}\big)|\geq \psi_n\big\}\geq \mathrm{1}\big\{|\widetilde{\Phi}-\Phi\big(f_{0}\big)|< \psi_n\big\}
$$
and, therefore,
\begin{equation}
\label{eq:ol13}
I_2\geq \bE^M\bE_{f_{y}}\Big[\mathrm{1}_{\cY}(y)\mathrm{1}\big\{|\widetilde{\Phi}-\Phi\big(f_{0}\big)|< \psi_n\big\}\Big].
\end{equation}
Next, recall that $\bE_f$ means the integration over $\cX^n=\cX\times\cdots\times\cX$ w.r.t. the product measure $\nu_n$ having the density
$$
\mathfrak{p}_f(x)=\prod_{i=1}^nf(x_i),\quad x=(x_1,\ldots,x_n)\in\cX^n.
$$
Recall also that an estimator  is a measurable function of observation. Hence, there exists a measurable function $\widetilde{G}:\cX^n\to\bR$ such that
$\widetilde{\Phi}=\widetilde{G}(X^{(n)})$.
Thus, for any $y\in\bY_\rho^*$ one has
$$
\bE_{f_{y}}\Big[\mathrm{1}\big\{|\widetilde{\Phi}-\Phi\big(f_{0}\big)|< \psi_n\big\}\Big]
=\int_{\cX^n}\Big[\mathrm{1}\big\{|\widetilde{G}(x)-\Phi\big(f_{0}\big)|< \psi_n\big\}\Big]\mathfrak{p}_{f_{y}}(x)\nu_n(\rd x).
$$
Moreover, in view of Lemma \ref{lem:bound_for_f_y} for any $y\in\bY^*_\rho$
$$
\mathfrak{p}_{f_{y}}(x)\geq e^{-1}\mathfrak{p}^*_y(x)=:\prod_{i=1}^nf^*_y(x_i),\quad\forall x\in\cX^n.
$$
Note that  $\mathfrak{p}^*_y(x)\geq 0$ for any $x\in\cX^n$ and $y\in\bY$, see Remark \ref{rem:after-lem:bound_for_f_y}. Therefore,
\begin{align*}
\mathrm{1}_{\cY}(y)\bE_{f_{y}}\Big[\mathrm{1}\big\{|\widetilde{\Phi}-\Phi\big(f_{0}\big)|&< \psi_n\big\}\Big]
\\
&\geq e^{-1}\mathrm{1}_{\cY}(y)\int_{\cX^n}\Big[\mathrm{1}\big\{|\widetilde{G}(x)-\Phi\big(f_{0}\big)|< \psi_n\big\}\Big]\mathfrak{p}^*_{y}(x)\nu_n(\rd x).
\end{align*}
Introduce the following quantities: for any $y\in\bY^M$
$$
Z=\bE^M(Z_y),\quad\;\; Z_y=\frac{\mathfrak{p}^*_{y}(X^{(n)})}{\mathfrak{p}_{f_0}(X^{(n)})},\quad\;\; \cI_y=\int_{\cX^n}\mathfrak{p}^*_{y}(x)\nu_n(\rd x)
=\bigg[\int_{\cX}f^*_y(w)\nu(\rd w)\bigg]^n.
$$
It is worth noting that $Z_y$ is well-defined in view of Assumption \ref{ass:density_0} and that $\cI_y\geq 0$ for any $y\in\bY$.

\noindent Thus, denoting $\overline{\cY}=\bY^M\setminus\cY$ we get
\begin{equation*}
\mathrm{1}_{\cY}(y)\bE_{f_{y}}\Big[\mathrm{1}\big\{|\widetilde{\Phi}-\Phi\big(f_{0}\big)|< \psi_n\big\}\Big]
\geq e^{-1}\bE_{f_{0}}\Big[Z_y\mathrm{1}\big\{|\widetilde{\Phi}-\Phi\big(f_{0}\big)|< \psi_n\big\}\Big]-e^{-1}\mathrm{1}_{\overline{\cY}}(y)\cI_y.
\end{equation*}
It yields together with (\ref{eq:ol13})
\begin{equation*}
eI_2\geq \bE_{f_{0}}\Big[Z\;\mathrm{1}\big\{|\widetilde{\Phi}-\Phi\big(f_{0}\big)|< \psi_n\big\}\Big]-\bE^M\big\{\mathrm{1}_{\overline{\cY}}(y)\cI_y\big\}
\end{equation*}
and applying the Cauchy-Schwarz inequality we get
\begin{equation}
\label{eq:ol14}
eI_2\geq \bE_{f_{0}}\Big[Z\;\mathrm{1}\big\{|\widetilde{\Phi}-\Phi\big(f_{0}\big)|< \psi_n\big\}\Big]-\sqrt{\bP^M\big(\overline{\cY}\big)\bE^M\big\{\cI^2_y\big\}}.
\end{equation}
At last, in view of Assumptions \ref{ass:density_xi_in_class}, \ref{ass:on-Phi} and the Tchebychev inequality, one has
$$
\bP^M\big(\overline{\cY}\big)\leq \bP^M\big(\overline{\bY^*_\rho}\big)+\bP^M\big(\overline{\bY_{\cF}}\big)+\bP^M\big(\overline{\bY_\Phi(\psi_n)}\big)\leq \tfrac{1}{144}+\tfrac{2}{231}\leq \tfrac{1}{64}.
$$
We deduce finally from (\ref{eq:ol14}), putting for brevity $R=8^{-1}\sqrt{\bE^M\big\{\cI^2_y\big\}}$
\begin{equation}
\label{eq:ol15}
eI_2\geq \bE_{f_{0}}\Big[Z\;\mathrm{1}\big\{|\widetilde{\Phi}-\Phi\big(f_{0}\big)|< \psi_n\big\}\Big]-R.
\end{equation}

$\mathbf{1^0c}.\;$ We deduce from (\ref{eq:ol12}) and (\ref{eq:ol15}) that for any estimator $\widetilde{\Phi}$
\begin{eqnarray*}
eR_n\big(\widetilde{\Phi}\big)&\geq& eI_1+eI_2\geq I_1+eI_2
\nonumber\\*[2mm]
&\geq& \alpha^2_n(c)\bE_{f_{0}}\Big[\mathrm{1}\big\{|\widetilde{\Phi}-\Phi\big(f_{0}\big)|\geq \psi_n\big\}\Big]+\bE_{f_{0}}\Big[Z\;\mathrm{1}\big\{|\widetilde{\Phi}-\Phi\big(f_{0}\big)|< \psi_n\big\}\Big]-R
\nonumber\\*[2mm]
&\geq&\bE_{f_{0}}\big[Z\wedge \alpha^2_n(c)\big]-R.
\end{eqnarray*}
Since the right hand side of obtained inequality is independent of $\widetilde{\Phi}$ we get
\begin{eqnarray}
\label{eq:ol16}
R_n:=\inf_{\widetilde{\Phi}}R_n\big(\widetilde{\Phi}\big)\geq e^{-1}\Big[\bE_{f_{0}}\big[Z\wedge \alpha^2_n(c)\big]-R\Big].
\end{eqnarray}

$\mathbf{2^0}.\;$ We have
\begin{equation*}
\mI:=\int_{\cX}f^*_y(w)\nu(\rd w)=\big[1-\Sigma_M\big]e^{-\frac{\rho_y-\Sigma_M}{1-\Sigma_M}}+\rho_y.
\end{equation*}
Hence, if $\lambda_m=0$ for all $m=1,\ldots,M$ then $\int_{\cX}f^*_y(w)\nu(\rd w)=1$ for all $y\in\bY^M$ and, therefore,
\begin{equation}
\label{eq:ol17}
R=\tfrac{1}{8}.
\end{equation}
Assume now that $\bY\subseteq[0,1]^M$. Since $\rho_y\geq 0$ and $\Sigma_M\leq 0.25$ in view of Assumption \ref{ass:one-of-two}  we have
$$
\tfrac{\rho_y-\Sigma_M}{1-\Sigma_M}\geq -1/3.
$$
Noting that   $e^{-t}\leq 1-t+\frac{1}{2}t^2$ for all $t\geq -1$,
 we have
$$
e^{-\frac{\rho_y-\Sigma_M}{1-\Sigma_M}}\leq 1-\tfrac{\rho_y-\Sigma_M}{1-\Sigma_M}+\tfrac{(\rho_y-\Sigma_M)^2}{2(1-\Sigma_M)^2}.
$$
It yields
\begin{equation*}
\mI\leq 1-\Sigma_M-(\rho_y-\Sigma_M)+\tfrac{(\rho_y-\Sigma_M)^2}{2(1-\Sigma_M)}+\rho_y\leq 1+\tfrac{2}{3}(\rho_y-\Sigma_M)^2.
\end{equation*}
Hence
$$
\bE^M\big\{\cI^2_y\big\}\leq \bE^M\Big\{\big[1+\tfrac{2}{3}(\rho_y-\Sigma_M)^2\big]^{2n}\Big\}\leq 2.
$$
The last inequality follows from Lemma \ref{lem:probab} applied with $K=M$, $J=2n$, $b=2/3$, $T=1$ and $\chi_K=\rho_y-\Sigma_M$. The assumption of the lemma follows from the condition $256\mS_M^2Mn\leq 1$ given in (\ref{eq:ass-on-parameters}) of Assumption \ref{ass:one-of-two}. It yields
\begin{equation}
\label{eq:ol18}
R\leq 8^{-1}\sqrt{2}<\tfrac{1}{4}.
\end{equation}
We obtain from (\ref{eq:ol16}), (\ref{eq:ol17}) and (\ref{eq:ol18})
\begin{eqnarray}
\label{eq:ol19}
R_n\geq e^{-1}\Big[\bE_{f_{0}}\big[Z\wedge \alpha^2_n(c)\big]-\tfrac{1}{4}\Big].
\end{eqnarray}

$\mathbf{3^0}.\;$ We obviously have
$$
2[Z\wedge \alpha^2_n(c)]=Z+\alpha^2_n(c)-\big|Z-\alpha^2_n(c)\big|.
$$
In view of the Fubini theorem
$$
\kappa:=\bE_{f_{0}}(Z):=\bE_{f_{0}}\big(\bE^M[Z_y]\big)=\bE^M\big(\bE_{f_{0}}[Z_y]\big)=:\bE^M[\cI_y].
$$
Noting that in view of Lemma \ref{lem:bound_for_f_y} $f_{y}(x)\leq f^*_y(x)$ for all $x\in\cX$ if  $y\in\bY^*_\rho$. It yields
$$
\int_{\cX}f^*_y(x)\nu(\rd x)\geq \int_{\cX}f_{y}(x)\nu(\rd x)=1,\quad \forall y\in\bY^*_\rho
$$
in view of (\ref{eq:int=1}). Since $\cI_y\geq 0$ for any $y\in\bY$ we obtain
\begin{equation*}
\kappa=\bE^M\big[\cI_y\big]\geq \bE^M\big[\mathrm{1}_{\bY^*_\rho}(y)\cI_y\big]\geq \bP^M(\bY^*_\rho)\geq \tfrac{143}{144}.
\end{equation*}
We get using Cauchy-Schwartz inequality
\begin{eqnarray}
\label{eq:ol20}
2\bE_{f_{0}}\big[Z\wedge \alpha^2_n(c)\big]&=&\kappa+\alpha^2_n(c)-\bE_{f_{0}}\big[|Z-\alpha^2_n(c)|\big]
\nonumber\\*[2mm]
&\geq& \kappa+\alpha^2_n(c)-
\sqrt{\bE_{f_{0}}\big[\big(Z-\alpha^2_n(c)\big)^2\big]}
\nonumber\\*[2mm]
&=&\kappa+\alpha^2_n(c)-
\sqrt{\bE_{f_{0}}\big[Z^2\big]-2\alpha^2_n(c)\kappa+\alpha^4_n(c)}
\nonumber\\*[2mm]
&=&\frac{\kappa^2+4\kappa\alpha^2_n(c)-\bE_{f_{0}}\big[Z^2\big]}{\kappa+\alpha^2_n(c)+
\sqrt{\bE_{f_{0}}\big[Z^2\big]-2\alpha^2_n(c)\kappa+\alpha^4_n(c)}}
\nonumber\\*[2mm]
&\geq&\frac{4\kappa-\alpha^{-2}_n(c)\bE_{f_{0}}\big[Z^2\big]}{1+\alpha^{-2}_n(c)+
\sqrt{\alpha^{-4}_n(c)\bE_{f_{0}}\big[Z^2\big]+1}}.
\end{eqnarray}
Thus, if we prove that
\begin{equation}
\label{eq:ol21}
\lim_{n\to\infty}\alpha^{-2}_n(c)\bE_{f_{0}}\big[Z^2\big]=0,
\end{equation}
then, since $\alpha_n(c)\to\infty$, we deduce from (\ref{eq:ol20}) that
$$
\lim_{n\to\infty}\bE_{f_{0}}\big[Z\wedge \alpha^2_n(c)\big]\geq\kappa.
$$
It yields together with (\ref{eq:ol19})
$$
\lim_{n\to\infty}R_n\geq \tfrac{107}{144e}
$$
that completes the proof of theorem.
\vskip0.1cm
$\mathbf{4^0}.\;$ It remains to establish (\ref{eq:ol21}). Note that
$$
Z^2=\big[\bE^M(Z_y)\big]^2=\bigg[\int_{\bR^M}Z_y\bP^M(\rd y)\bigg]^2=\int_{\bR^M}\int_{\bR^M}Z_y Z_z\bP^M(\rd y)\bP^M(\rd z).
$$
Hence in view of the Fubini theorem
\begin{eqnarray}
\label{eq:ol22}
\bE_{f_{0}}\big[Z^2\big]&=&\int_{\bR^M}\int_{\bR^M}\bE_{f_{0}}\big[ZyZ_z\big]\bP^M(\rd y)\bP^M(\rd z)
\nonumber\\*[2mm]
&=&\int_{\bR^M}\int_{\bR^M}\bigg[\int_{\cX}\tfrac{f^*_y(x)f^*_z(x)}{f_0(x)}\nu(\rd x)\bigg]^n\bP^M(\rd y)\bP^M(\rd z)
\nonumber\\*[2mm]
&
=:&\int_{\bR^M}\int_{\bR^M}\big[\cE(y,z)\big]^n\bP^M(\rd y)\bP^M(\rd z).
\end{eqnarray}
It is worth mentioning that $\cE(\cdot,\cdot)\geq 0$ in view of Assumption \ref{ass:one-of-two}.
To make easier  the rest of the proof we consider separately two cases.
\vskip0.1cm
$\mathbf{4^0a}.\;$ Let $\lambda_m=0$ for all $m=1,\ldots, M$. In this case $f^*_y\equiv f_{y}$ and we easily check that
$$
\cE(y,z)=1+\sum_{m=1}^M S_m y_m z_m, \quad \forall y,z\in\bR^M\quad\Rightarrow\quad \big[\cE(y,z)\big]^n\leq \prod_{m=1}^Me^{nS_m y_mz_m}.
$$
Since $\bP^M$ is  product measure we get from (\ref{eq:ol22})
$$
\bE_{f_{0}}\big[Z^2\big]\leq \prod_{m=1}^M\int_{\bR^2}e^{nS_m uv}\bP(\rd u)\bP(\rd v)
$$
and (\ref{eq:ol21}) follows from the condition (\ref{eq:main-assump}) of the theorem with $\mathbf{i}=1$.

\vskip0.1cm

$\mathbf{4^0b}.\;$ The computations in the remaining case are a  bit more cumbersome. Recall that $y\in[0,1]$ and $\rho_y\geq 0$. First we remark that in view of Assumptions \ref{ass:func_lambda} and \ref{ass:density_0}

\begin{eqnarray}
\label{eq:ol23}
\cE(y,z)&=&\big[1-\Sigma_M\big]^2e^{-\frac{\rho_y-\Sigma_M}{1-\Sigma_M}}e^{-\frac{\rho_z-\Sigma_M}{1-\Sigma_M}}
+ \sum_{m=1}^M S_m y_m z_m
\nonumber\\*[2mm]
&&+(1-\Sigma_M)\rho_ye^{-\frac{\rho_z-\Sigma_M}{1-\Sigma_M}}+(1-\Sigma_M)\rho_ze^{-\frac{\rho_y-\Sigma_M}{1-\Sigma_M}}
\nonumber\\*[2mm]
&=&\big[1-\Sigma_M\big]^2e^{-\frac{\Upsilon_y+\Upsilon_z}{1-\Sigma_M}}+ \sum_{m=1}^M S_m y_m z_m
\nonumber\\*[2mm]
&&+(1-\Sigma_M)\rho_ye^{-\frac{\Upsilon_z}{1-\Sigma_M}}+(1-\Sigma_M)\rho_ze^{-\frac{\Upsilon_y}{1-\Sigma_M}}
\nonumber\\*[2mm]
&\leq&\big[1-\Sigma_M\big]^2e^{-\frac{\Upsilon_y+\Upsilon_z}{1-\Sigma_M}}+ \sum_{m=1}^M S_m y_m z_m +2(1-\Sigma_M)\Sigma_M
\nonumber\\*[2mm]
&&+(1-\Sigma_M)\Upsilon_ye^{-\frac{\Upsilon_z}{1-\Sigma_M}}+(1-\Sigma_M)\Upsilon_ze^{-\frac{\Upsilon_y}{1-\Sigma_M}}.
\end{eqnarray}
Here we have used that $\rho_y$ is a linear function of $y$ and put  $\Upsilon_y=\rho_{y}-\Sigma_M$,
$\Upsilon_z=\rho_{z}-\Sigma_M$.

Since $\rho_{y}\geq 0$ and $\Sigma_M\leq 0.25$ in view of Assumption \ref{ass:one-of-two}  we have
$$
\tfrac{\Upsilon_y+\Upsilon_z}{1-\Sigma_M}\geq -2/3,\quad \tfrac{\Upsilon_y}{1-\Sigma_M}\geq -1/3,\quad \tfrac{\Upsilon_z}{1-\Sigma_M}\geq -1/3.
$$
Using once again  $e^{-t}\leq 1-t+\frac{1}{2}t^2$ for all $t\geq -1$, 
 we have
$$
e^{-\frac{\Upsilon_y+\Upsilon_z}{1-\Sigma_M}}\leq 1-\tfrac{\Upsilon_y+\Upsilon_z}{1-\Sigma_M}+\tfrac{(\Upsilon_y+\Upsilon_z)^2}{2(1-\Sigma_M)^2}.
$$
On other hand $|e^{-t}-1|\leq e|t|$ for all $t\geq -1$ and, therefore,
$$
\Big|e^{-\frac{\Upsilon_y}{1-\Sigma_M}}-1\Big|\leq \tfrac{e|\Upsilon_y|}{1-\Sigma_M},\quad\;\Big|e^{-\frac{\Upsilon_z}{1-\Sigma_M}}-1\Big|\leq \tfrac{e|\Upsilon_z|}{1-\Sigma_M}.
$$
We deduce from (\ref{eq:ol23})
\begin{eqnarray}
\cE(y,z)&\leq& 1-\Sigma_M^2+ \sum_{m=1}S_m y_m z_m+ 2^{-1}(\Upsilon_y+\Upsilon_z)^2+2e|\Upsilon_y||\Upsilon_z|
\nonumber\\*[2mm]
&\leq& 1+ \sum_{m=1}S_m y_m z_m+ (e+1)\Upsilon_y^2+(e+1)\Upsilon_z^2.
\end{eqnarray}
Using the trivial inequality
$$
(1+u+v+w)^n\leq 2^{-1}(1+2u)^n +4^{-1}(1+4v)^n+4^{-1}(1+4w)^n
$$
we have
\begin{eqnarray*}
\cE(y,z)&\leq& \tfrac{1}{2}\bigg(1+2\sum_{m=1}S_m y_m z_m\bigg)^n+\tfrac{1}{4}\Big(1+4(e+1)\Upsilon_y^2\Big)^n+\tfrac{1}{4}\Big(1+4(e+1)\Upsilon_z^2\Big)^n
\nonumber\\*[2mm]
&\leq &2^{-1}\prod_{m=1}^Me^{2n S_my_mz_m}+4^{-1}\Big(1+4(e+1)\Upsilon_y^2\Big)^n+4^{-1}\Big(1+4(e+1)\Upsilon_z^2\Big)^n.
\end{eqnarray*}
The inequality (\ref{eq:ol21}) follows now from the condition (\ref{eq:main-assump}) of the theorem with $\mathbf{i}=2$ and Lemma \ref{lem:probab}
applied with $J=n$, $K=M$, $T=1$, $b=4(e+1)$ and $\chi_K=\Upsilon_y,\Upsilon_z$ under  condition $256\mS_M^2Mn\leq 1$ of Assumption \ref{ass:one-of-two}.
The theorem is proved.
\epr

\subsection{Proof of Theorem \ref{th:adaptive-lb-bounded}}
We will apply Theorems \ref{th:adaptive-rate} and \ref{th:adaptive-rate_K}.

\vskip0.1cm

$\mathbf{1^0}.\;$ Let us  prove that $\Psi^b$ is an adaptive rate of convergence on $\Theta^b_{(r),1}$. First of all we remark that establishing Proposition \ref{prop:nonexitence-1} we already proved that all assumptions of Theorem \ref{th:generic-adaptive} are checked on the parameter set including $\Theta^b_{(r),1}$. Thus, in order to apply Theorem \ref{th:adaptive-rate} it suffices to verify Conditions $\mathbf{A_1}-\mathbf{A_3}$.

\vskip0.1cm

$\mathbf{1^0a}.\;$ Verification of Condition $\mathbf{A_1}$. Fix $\theta=(\blb,\bbr)$ and denote for brevity
$
U=\frac{\blb\bbr}{\blb\bbr+d(\bbr-1)}.
$
Then, the set $\Theta^0[\theta]$ consists of couples $\theta^\prime=(\bga,\bbs)$ verifying
$$
\frac{\bga\bbs}{\bga\bbs+d(\bbs-1)}=U\quad\Leftrightarrow\quad \bga=\tfrac{dU(1-1/\bbs)}{1-U}.
$$
Thus, $\bbs\mapsto \bga(\bbs)$ is a continuous function of $\bbs$ on $[2,\infty]$ and obviously $\Theta^0[\theta]\in\mT^{1}(\Theta^b_{(r),1})$ whatever $b>0$.
On the other hand the set $\Theta^+[\theta]$ consists of couples $\theta^\prime=(\bga,\bbs)$ verifying
$$
\tfrac{1}{2}>\frac{\bga\bbs}{\bga\bbs+d(\bbs-1)}> U
$$
and this set obviously contains a ball of $\bR^2$. Condition $\mathbf{A_1}$ is checked.

\vskip0.1cm

$\mathbf{1^0b}.\;$ Verification of Condition $\mathbf{A_2}$. Since
$$
\varphi_n(\theta)=\big(\tfrac{1}{n}\big)^{\mz^*(\theta)},\quad\;\varphi_n(\theta^\prime)=\big(\tfrac{1}{n}\big)^{\mz^*(\theta^\prime)}
$$
and $\mz^*(\theta^\prime)>\mz^*(\theta)$, because $\theta^\prime\in\Theta^+[\theta]$,
choosing $c(\theta,\theta^\prime)=\frac{1}{2}+\frac{\mz^*(\theta)}{2\mz^*(\theta^\prime)}<1$ we obtain
$$
\alpha_n\big(\theta,\theta^\prime\big)=\big(\tfrac{1}{n}\big)^{\frac{1}{2}[\mz^*(\theta)-\mz^*(\theta^\prime)]}\to\infty,\quad n\to\infty.
$$
 Condition $\mathbf{A_2}$ is checked.

 \vskip0.1cm

 $\mathbf{1^0c}.\;$ Verification of Condition $\mathbf{A_3}$. Since $\mz^*(\theta)<1/2$ for any $\theta\in\Theta^b_{(r),1}$ we have
 $$
\limsup_{n\to\infty}\mathfrak{p}^2_n\big[\varphi_n(\theta^\prime)\big]^{1-c(\theta,\theta^\prime)}\leq \limsup_{n\to\infty} \sqrt{\ln(n)}\big(\tfrac{1}{n}\big)^{\frac{1}{2}[\mz^*(\theta^\prime)-\mz^*(\theta)]}=0.
$$
 Condition $\mathbf{A_3}$ is checked. Since we assumed that $\Psi^b$ is admissible on $\Theta^b_{(r),1}\cup\Theta^b_{(r),2}$ and, therefore on $\Theta^b_{(r),1}$ applying Theorem \ref{th:adaptive-rate} we assert that $\Psi^b$ is an adaptive rate of convergence on $\Theta^b_{(r),1}$.

\vskip0.1cm

$\mathbf{2^0}.\;$ Let us  prove that $\Psi^b$ is an adaptive rate of convergence on $\Theta^b_{(r),2}$. First of all we remark that establishing Proposition \ref{prop:nonexitence-2} we already proved that all assumptions of Theorem \ref{th:generic-adaptive} are checked on the parameter set including $\Theta^b_{(r),2}$. Indeed, as we already mentioned, see Remark \ref{rem5}, in the isotropic case the requirement $\rho_{\theta,\theta^\prime}\geq 1$ is automatically fulfilled. Also, the set on which  Proposition \ref{prop:nonexitence-2} is proved contains the requirement $q^{\prime}\leq q$ which is fulfilled in the considered case since $q^{\prime}=q=\infty$. The verification of Conditions $\mathbf{A_1}-\mathbf{A_3}$ can be done following the same lines as in $\mathbf{1^0a}-\mathbf{1^0c}$. Applying Theorem \ref{th:adaptive-rate} we assert that $\Psi^b$ is an adaptive rate of convergence on $\Theta^b_{(r),2}$.

The assertion of the theorem follows now from Theorem \ref{th:adaptive-rate_K} with $K=2$.
\epr

\subsection{Proofs of Propositions \ref{prop:nonexitence-1} and \ref{prop:nonexitence-2}.}

We start this section with the choice of elements used in the construction (\ref{eq:family}) and  with formulation some auxiliary results.
Furthermore, $c_1,c_2,\ldots, C_1,C_2,\ldots$ stands for constant independent of $n$.

\subsubsection{Preliminary remarks.}
\label{sec:sunsec-preliminary_remarks}

\paragraph*{Parameterized family of functions.}
Fix vector $\vec{\sigma}=(\sigma_1,\ldots,\sigma_d)\in (0,1]^d$ and set for any $x=(x_1,\ldots,x_d)\in\bR^d$
$$
\Pi_x=\big(x_{1}-\sigma_1/2,x_{1}+\sigma_1/2\big)\times\cdots\times\big(x_{d}-\sigma_d/2,x_{d}+\sigma_d/2\big).
$$
In words, $\Pi_x$ is a rectangle in $\bR^d$ centered at $x$ with edges of lengths $\sigma_1, \ldots \sigma_d$,
that are parallel to the coordinate axes.
Fix $M\in\bN^*$ and let $\{x_1,\ldots,x_M\}\subset\bR^d$ be such that
\begin{equation}
\label{eq:disjoint-Pi}
\Pi_{x_m}\cap \Pi_{x_l}=\emptyset, \quad \forall l,m=1,\ldots,M.
\end{equation}
Let $\Lambda:\bR^d\to \bR$,  $\Lambda\in\bC^{\infty}(\bR^d)$  be a function satisfying the following conditions:
\begin{equation}
\label{eq:Lambda}
\Lambda(x)=0,\; \; \forall x\notin [-1/2,1/2]^d,\quad \int_{\bR^d}\Lambda(x)\rd x =0,\quad \|\Lambda\|_\infty\leq 1.
\end{equation}
Let $A>0$ be chosen later.
Define for any $m=1,\ldots, M$
\begin{equation}
\label{eq:Pi-m}
\Lambda_m(x)=A\Lambda\big((x-x_{m})/\vec{\sigma}\big),\quad \Pi_{x_m}=
\Big\{x\in\bR^d:\; \big|(x-x_{m})/\vec{\sigma}\big|_\infty\leq 1\Big\},
\end{equation}
where the division is understood in the coordinate--wise sense.
It is obvious that $\Lambda_m$ is supported on $\Pi_{x_m}$ for any  $m=1,\ldots, M$, and by construction  $\Pi_{x_m}$ are disjoint.

\vskip0.1cm

Let $f_0$ be a probability density on $\bR^d$. Choose $\bY=\{-1,1\}$ and, therefore,   $\bY^M=\{-1,1\}^M$.
\begin{eqnarray*}
&&
f_y(x)=f_0(x)+A\sum_{m=1}^M y_m\Lambda\big((x-x_{m})/\vec{\sigma}\big),\quad x\in\bR^d.
\end{eqnarray*}
Thus, we come to the collection (\ref{eq:family}) with $\rho_y=0$ since in view of (\ref{eq:Lambda}) $\lambda_m=0$ for all $m=1,\ldots,M$.
We constat also that Assumption \ref{ass:func_lambda} is verified.

Set $\bma=\prod_{l=1}\sigma_l$. The following result was proved in \cite{gl14}.

\begin{lemma}
\label{lem:lem-from-GL14}
 For any $M\in\bN^*$,   $y\in \{-1,1\}^M$, $\vec{\beta},\vec{r}$, $\vec{L}$, $\vec{\sigma}$ and $M$ the function $x\mapsto A\sum_{m=1}^M y_m\Lambda\big((x-x_{m})/\vec{\sigma}\big)$ belongs to $\bN_{\vec{r},d}(\vec{\beta},\frac{1}{2}\vec{L})$ if for some $C_1>0$, completely determined by $\Lambda$, $\vec{\beta}$ and $\vec{r}$
\begin{eqnarray}
\label{eq:cond-lemma-fromGL14}
&&A\sigma_l^{-\beta_l}
\big[\bma M\big]^{1/r_l}\leq C_1  L_l,\quad\forall l=1,\ldots, d.
\end{eqnarray}

\end{lemma}
\noindent Also, simplest algebra shows that for any $M\in\bN^*$ and  $y\in \{-1,1\}^M$ the function $x\mapsto A\sum_{m=1}^M y_m\Lambda\big((x-x_{m})/\vec{\sigma}\big)$ belongs to $\bB_q(Q/2)$ if
for some $C_2>0$, completely determined by $\Lambda$,
\begin{equation}
\label{eq:cond_q-norm}
A(M\bma)^{1/q}\leq C_2.
\end{equation}

\paragraph*{First construction of $f_0$.} Let $R:=\int_{-1}^1e^{-\frac{1}{1-z^2}}\rd z$, and
$$
U(x):=R^{-d}e^{-\sum_{j=1}^d\frac{1}{1-x_i^2}}\mathrm{1}_{[-1,1]^d}(x),\;\; x=(x_1,\ldots,x_d)\in\bR^d.
$$
For  $N\geq 1$ and $a>0$ define
\begin{gather*}
\bar{f}_{0,N}(x):=N^{-d}\int_{\bR^d}U(y-x)\mathrm{1}_{[1,N+1]^d}(y)\rd y,\;\;\;
f_{0,N,a}(x):=a^d\bar{f}_{0,N}\big(xa\big).
\end{gather*}
We will need the following result.
\begin{lemma}[\cite{GL20a}, Lemma 1]
\label{lem:density-f_0,N}
For any $N\geq 1$ and $a>0$, $f_{0,N,a}$ is a probability density.

For any $\vec{\beta}, \vec{L}\in (0,\infty)^d$ and  $\vec{r}\in [1,\infty]^d$
there exists $a_0>0$ such that for any $a\leq a_0$
$$
f_{0,N,a}\in\bN_{\vec{r},d}\big(\vec{\beta},\tfrac{1}{2}\vec{L}\big),\quad\forall N\geq 1.
$$
\end{lemma}

\noindent Applying the Young inequality, see \cite{Folland}, Theorem 8.9., we have for any $a>0$,  $N\geq 1$ and $u> 1$
\begin{equation}
\label{eq:f_0-in-L_q}
\|f_{0,N,a}\|_u\leq a^{d(1-1/u)}N^{-d+d/u}\leq N^{-d(1-1/u)},\quad \forall a\leq 1.
\end{equation}
Let $\vartheta=(\vec{\beta},\vec{r},q,\vec{L},Q)=:(\theta,\bar{\theta})\in\Theta_{(r)}\times\Theta_{(c)}$ and $\vartheta^\prime=(\vec{\gamma},\vec{s},q^\prime,\vec{L}^{\prime},Q^\prime)=:(\theta^\prime,\bar{\theta}^\prime)$ such that $\theta^\prime\in\Theta[\theta]$
be given. Then in view of Lemma \ref{lem:lem-from-GL14} and (\ref{eq:f_0-in-L_q}) there exist $\mathbf{a}<1$ and $\mathbf{N}=\mathbf{N}(q,q^\prime,Q,Q^\prime)>8$ such that for any $N\geq \mathbf{N}$
\begin{equation}
\label{eq:family_0}
f_{0,N,\mathbf{a}}\in\bN_{\vec{r},d}\big(\vec{\beta},\tfrac{1}{2}\vec{L}\big)\cap\bB_q\big(\tfrac{1}{2}Q\big)
\cap\bN_{\vec{s},d}\big(\vec{\gamma},\tfrac{1}{2}\vec{L}^\prime\big)\cap\bB_{q^\prime}\big(\tfrac{1}{2}Q^\prime\big)=:\bF_{\vartheta,\vartheta^\prime}.
\end{equation}
Set $f_{0,N}(x)=f_{0,N,\mathbf{a}}(x)$ and choose $f_0=f_{0,N}$. Thus, we come to the collection
\begin{equation}
\label{eq:family_1}
f_y(x)=f_{0,N}(x)+A\sum_{m=1}^M y_m\Lambda\big((x-x_{m})/\vec{\sigma}\big),\quad x\in\bR^d, \quad y\in\{-1,1\}^M, \; N\geq \mathbf{N}.
\end{equation}
This collection of functions will be used in the proof of Proposition \ref{prop:nonexitence-1}, where the parameters $N,A$ and $\vec{\sigma}$ will be chosen. As we already mentioned
\begin{equation}
\label{eq:f_0-in-2classes-1}
f_{0,N}\in \bF_{\vartheta,\vartheta^\prime},\quad \forall N\geq \mathbf{N}.
\end{equation}

\paragraph*{Second construction of $f_0$.} Set $\varsigma^\prime=(\vec{\gamma},\vec{s})$ and recall that $\varsigma=(\vec{\beta},\vec{r})$. Define for any $q^\prime\leq q$
$$
\rho_l=\Big[\tfrac{1}{\beta_l}\big(\tfrac{1}{r_l}-\tfrac{1}{q}\big)\Big]\wedge\Big[\tfrac{1}{\gamma_l}\big(\tfrac{1}{s_l}-\tfrac{1}{q^\prime}\big)\Big],
$$
and recall that for any $\theta^\prime=(\varsigma^\prime,q^\prime)\in\Theta^{\prime\prime}$
\begin{equation}
\label{eq:assump-rho}
\rho_{\theta,\theta^\prime}:=\sum_{l=1}^d\rho_l\geq 1.
\end{equation}
Since, furthermore $\theta^\prime,\theta$ are fixed, we will write for brevity $\rho$ instead of $\rho_{\theta,\theta^\prime}$.

\noindent For any  $t\in(0,1]$ set $\vec{\mh}=(\mh_1,\ldots,\mh_d)$,
$
\mh_l=t^{\frac{\rho_l}{\rho}}
, l=1,\ldots,d,
$
and consider the following function
$$
G_t(x)=\mathbf{c}t^{-1/q}f_{0,\mathbf{N},\mathbf{a}}\big(x/\vec{\mh}\big),\quad x\in\bR^d.
$$
Here $0<\mathbf{c}<1$ is sufficiently small constant independent on $n$ and $t$ whose choice will done below.
We remark that in view of Lemma \ref{lem:density-f_0,N} the assertion of Lemma \ref{lem:lem-from-GL14} is applicable to $G_t$ with $A=\mathbf{c}t^{-1/q}=\mathbf{c}\big(\prod_{l=1}^d\mh_l\big)^{-1/q}$, $M=1$ and $\sigma_l=\mh_l$, $l=1,\ldots d$. Here we have taken into account that $\vec{\mh}\in (0,1]^d$ since $t\in(0,1]$.
Hence this function belongs to $\bN_{\vec{s},d}\big(\vec{\gamma},\tfrac{1}{2}\vec{L}^\prime\big)$ if for some constant $C^\prime_1$ completely determined by $\vec{\gamma}$, $\vec{s}$, $\mathbf{N}$ and $\mathbf{a}$ the following inequalities hold.
$$
\mathbf{c}\mh_l^{-\gamma_l}t^{1/s_l-1/q}\leq C_1L^\prime_l, \quad l=1,\ldots,d.
$$
Here we have used that $\prod_{l=1}^d\mh_l=t$. We assert that the latter inequality holds for any $q^\prime\leq q$ if
$$
\mathbf{c}t^{(1/s_l-1/q^\prime)(1-1/\rho)}\leq C_1L^\prime_l.
$$
Here we have also taken into account that $t\leq 1$.
Since $s_l\leq q^\prime$ and $\rho\geq 1$ in view of (\ref{eq:assump-rho}), choosing $\mathbf{c}\leq C_1\min_{l=1,\ldots, d}\{L_l\wedge L^\prime_l\}$ we assert that
$$
G_t\in\bN_{\vec{s},d}\big(\vec{\gamma},\tfrac{1}{2}\vec{L}^\prime\big), \quad \forall t\in (0,1).
$$
Note also that in view of (\ref{eq:f_0-in-L_q}) with $u=q^\prime$ we have $G_t\in\bB_{q^\prime}\big(\tfrac{1}{2}Q^\prime\big)$ for any $q^\prime\leq q$.

Since the definition of $\rho_l, l=1,\ldots$ is "symmetric" in $(\vec{\gamma},\vec{s}, q^{\prime})$ and  $(\vec{\beta},\vec{r}, q)$ we can assert finally, taking into account  the choice of constant $\mathbf{c}$,  that
\begin{equation}
\label{eq:family_0-second}
G_t\in\bN_{\vec{r},d}\big(\vec{\beta},\tfrac{1}{2}\vec{L}\big)\cap\bB_q\big(\tfrac{Q}{2}\big)
\cap\bN_{\vec{s},d}\big(\vec{\gamma},\tfrac{1}{2}\vec{L}^\prime\big)\cap\bB_{q^\prime}\big(\tfrac{1}{2}Q^\prime\big)=\bF_{\vartheta,\vartheta^\prime},
\quad\forall t\in (0,1).
\end{equation}
Define, at last, for $t\in (0,1)$,
$$
f_{0,t}(x)=\big(1-\mathbf{c}t^{1-1/q}\big)f_{0,e^{-n},\mathbf{a}}(-x)+G_t(x), \quad x\in\bR^d.
$$
We will assume that $n$ is sufficiently large  in order to provide $e^{-n}<\mathbf{N}$. Obviously,  $f_{0,t}$ is a probability density for all $t\in(0,1]$
and in view of (\ref{eq:family_0}) and (\ref{eq:family_0-second}) we assert that for any $q^\prime\leq q$
\begin{equation}
\label{eq:f_0-in-2classes-3}
f_{0,t}\in \bN_{\vec{r}, d}(\vec{\beta}, \vec{L})\cap\bB_q(Q)\cap\bN_{\vec{s}, d}(\vec{\gamma}, \vec{L}^\prime)\cap\bB_q(Q^\prime),\quad
\forall t\in (0,1).
\end{equation}
Thus, we come to the following construction.
\begin{equation}
\label{eq:family_2}
f_y(x)=f_{0,t}(x)+A\sum_{m=1}^M y_m\Lambda\big((x-x_{m})/\vec{\sigma}\big),\quad x\in\bR^d, \quad y\in\{-1,1\}^M, \; t\in (0,1).
\end{equation}
This collection of functions will be used in the proof of Proposition \ref{prop:nonexitence-2}, where the parameters $t,A$ and $\vec{\sigma}$ will be chosen.

\subsubsection{Proof of Proposition \ref{prop:nonexitence-1}.}
We start with choosing the set $\cX=\{x_1,\ldots,x_M\}$. First we note that if  $\boldsymbol{\cX}$ will be the set with maximal cardinality, say $\boldsymbol{M}$, verifying (\ref{eq:disjoint-Pi})  and such that $\Pi_{x_m}\subset [2,N-1]^d$ for any $m=1,\ldots,M$, then
$$
(N-3)^{d}\geq \boldsymbol{M}\bma\geq (N-4)^d,\quad \forall N\geq 8,\;\;\forall \vec{\sigma}\in (0,1]^d
$$
It allows us to choose $M$ as follows. Let $0<\kappa<(1/2)^d$ be the number which will be chosen later. Set
\begin{equation}
\label{eq:connectionM&N}
 M\bma=\kappa N^d,
\end{equation}
and let $\Pi_{x_m}\subset [1,N]^d$ for any $m=1,\ldots,M$ and verify (\ref{eq:disjoint-Pi}).
Furthermore,  constants $c_1,c_2\ldots,$ below will be  independent of $\kappa$ as well.

\vskip0.1cm

$\mathbf{1^0}.\;$ In view of the choice of $\cX$ one has $\cup_{m=1}^M\Pi_{x_m}\subset[1,N]^d$ and, therefore,
$$
f_{0,N}(x)=\mathbf{a}^{d}N^{-d},\quad \forall x\in \cup_{m=1}^M\Pi_{x_m}.
$$
Firstly, we assert that Assumption \ref{ass:density_0} is fulfilled.

Secondly, it allows us to assert that for any $m=1,\ldots, M$
\begin{equation}
\label{eq1:proof-prop1}
S_m=\int_{\cX} \tfrac{\Lambda^2_m(x)}{f_0(x)}\nu(\rd x)=N^{d}\mathbf{a}^{-d}\int_{\bR^d} \Lambda^2_m(x)\rd x
=A^2N^{d}\mathbf{a}^{-d}\bma \|\Lambda\|_2^2= c_1A^2N^{d}\bma.
\end{equation}
Thirdly, choosing
\begin{equation}
\label{eq2:proof-prop1}
N^d=\mathbf{a}^{d}A^{-1}.
\end{equation}
and taking into account that $\|\Lambda\|_\infty\leq 1$ in view (\ref{eq:Lambda}), we can assert that
$
f_y
$
given in (\ref{eq:family_1}) satisfies
$$
f_y(x)\geq 0, \quad\forall x\in\bR^d,\;\;\forall y\in\{-1,1\}^M.
$$
Hence, Assumption \ref{ass:one-of-two} holds.

At last, putting for brevity $\kappa_1=\mathbf{a}^{d}\kappa$, we deduce from (\ref{eq:connectionM&N}) and (\ref{eq2:proof-prop1})
\begin{equation}
\label{eq200:proof-prop1}
 M\bma=\bka_1 A^{-1}.
\end{equation}

$\mathbf{2^0}.\;$  Choose $A=\left(n^{-1}\sqrt{\ln(n)}\right)^{2\mz(\theta)}$, where $\mz(\theta)$ is given in (\ref{eq:rate-exponent}), and set
\begin{eqnarray*}
\sigma_l&=&C_1^{-\frac{1}{\beta_l}}L_l^{-\frac{1}{\beta_l}}A^{\frac{1}{\beta_l}-\frac{1}{\beta_lr_l}}\kappa_1^{\frac{1}{\beta_lr_l}};
\\*[2mm]
M&=&C_1^{\frac{1}{\beta}}\BL A^{-\tau_{\varsigma}(1)}\kappa_1^{1-\frac{1}{\omega_\varsigma}},\quad\; \BL=\prod_{l=1}^d L_l^{\frac{1}{\beta_l}}.
\end{eqnarray*}
First we remark that (\ref{eq200:proof-prop1}) holds. Also we constat that
the condition
 (\ref{eq:cond-lemma-fromGL14}) of Lemma \ref{lem:lem-from-GL14} holds. It yields together with (\ref{eq:f_0-in-2classes-1})
 $$
 f_y\in\bN_{\vec{r},d}\big(\vec{\beta},\vec{L}\big)\cap\bB_q\big(Q\big),\quad \forall y\in\{-1,1\}^M.
 $$
 Since (\ref{eq:f_0-in-2classes-1}) also guarantees that $f_{0,N}\in\bN_{\vec{s},d}\big(\vec{\gamma},\vec{L}^\prime\big)\cap\bB_{q^\prime}\big(Q^\prime\big)$ we can assert that Assumption \ref{ass:density_xi_in_class} is verified.

 It remains to check that $\sigma_l\leq 1$ for any $l=1,\ldots,d$. Indeed, since $A<1$ for all $n$ large enough and $r_l\geq 1$ we have
 $$
 \sigma_l\leq C_1^{-\frac{1}{\beta_l}}L_l^{-\frac{1}{\beta_l}}\kappa_1^{\frac{1}{\beta_lr_l}}.
 $$
 Hence, choosing $\kappa$ sufficiently small, guaranteeing
 $
 \kappa_1\max_{l=1,\ldots,d}C_1^{-r_l}L_l^{-r_l}\leq 1
 $
 we come to the required constraints.

\smallskip

 $\mathbf{3^0}.\;$ Taking into account that $\int\Lambda=0$ in view (\ref{eq:Lambda}) we get
 $$
 \|f_y\|_2^2=\|f_{0,N}\|_2^2+A^2\bma M\|\Lambda\|_2^2,\quad \|f_{0,N}\|_2\leq  N^{-d/2}=c_2\sqrt{A},\quad \forall y\in\{-1,1\}^M.
 $$
 The last inequality follows from (\ref{eq:f_0-in-L_q}) with $u=2$ and the choice (\ref{eq2:proof-prop1}). It yields
 $$
\big|\|f_y\|_2-\|f_{0,N}\|_2\big|=\tfrac{\left|\|f_y\|^2_2-\|f_{0,N}\|^2_2\right|}{\|f_y\|_2+\|f_{0,N}\|_2}
\geq c_3\tfrac{A^2\bma M}{\sqrt{A}+A\sqrt{M\bma}}\quad \forall y\in\{-1,1\}^M.
 $$
Using (\ref{eq:connectionM&N}) we finally get
\begin{equation*}
\big|\|f_y\|_2-\|f_{0,N}\|_2\big|\geq c_4\kappa\sqrt{A}\quad \forall y\in\{-1,1\}^M.
\end{equation*}
We assert that Assumption \ref{ass:on-Phi} is fulfilled with
\begin{equation}
\label{eq3:proof-prop1}
\psi_n=2^{-1}c_{4}\kappa\Big(\tfrac{\sqrt{\ln(n)}}{n}\Big)^{\mz(\theta)}.
\end{equation}

 $\mathbf{4^0}.\;$ Choose the probability measure $P$ as follows:  $P\{1\}=P\{-1\}=1/2.$ It yields together with (\ref{eq1:proof-prop1})
 $$
 \int_{\bR^2}e^{\mathbf{i} nS_m uv}\bP(\rd u)\bP(\rd v)=\tfrac{1}{2}\big[e^{c_1\mathbf{i}nA^2N^{d}\bma}+e^{-c_1\mathbf{i}nA^2N^{d}\bma}\big]
 $$
 Note that in view of (\ref{eq:connectionM&N}) and (\ref{eq2:proof-prop1})
 $$
 A^2N^{d}\bma=\tfrac{A^2N^{d}M\bma}{M}=\tfrac{A^2N^{2d}\kappa}{M}=\tfrac{\mathbf{a}^{d}\kappa_1}{M}=\mathbf{a}^{d}C_1^{-\frac{1}{\beta}}\BL^{-1} A^{\tau_{\varsigma}(1)}\kappa_1^{\frac{1}{\omega_\varsigma}}=c_5A^{\tau_{\varsigma}(1)}\kappa_1^{\frac{1}{\omega_\varsigma}}.
 $$
It yields
$
nA^2N^{d}\bma=c_5\kappa_1^{1/\omega_\varsigma}\tfrac{\ln(n)}{n}
$
and, therefore,
$$
\int_{\bR^2}e^{\mathbf{i} nS_m uv}\bP(\rd u)\bP(\rd v)\leq 1+\tfrac{c_6\kappa_1^{2/\omega_\varsigma}\ln^2(n)}{n^2},
$$
for all $n$ large enough. It yields
\begin{equation*}
\prod_{m=1}^M \int_{\bR^2}e^{\mathbf{i} nS_m uv}\bP(\rd u)\bP(\rd v)\leq \Big( 1+\tfrac{c_6\kappa_1^{2/\omega_\varsigma}\ln^2(n)}{n^2}\Big)^M
\leq e^{ \tfrac{c_6M\kappa_1^{2/\omega_\varsigma}\ln^2(n)}{n^2}}.
\end{equation*}
It remains to note that
$$
\tfrac{c_6M\kappa_1^{2/\omega_\varsigma}\ln^2(n)}{n^2}=c_7\kappa_1^{1+1/\omega_\varsigma}\ln(n)
$$
and we get
\begin{equation}
 \label{eq4:proof-prop1}
\prod_{m=1}^M \int_{\bR^2}e^{\mathbf{i} nS_m uv}\bP(\rd u)\bP(\rd v)\leq n^{c_7\kappa_1^{1+1/\omega_\varsigma}}.
\end{equation}

 $\mathbf{5^0}.\;$  Let $\theta^\prime\in \Theta^{\prime}[\theta]$ and  $\mz(\theta)<\alpha< \mz(\theta^\prime)$ be fixed. Set $c=\frac{\alpha}{\mz(\theta^\prime)}$ and note that $c<1$. On the other hand
 $$
 \alpha_n(c):=\tfrac{\phi_n(\cF_\vartheta)}{\phi^{c}_n(\cF_{\vartheta^\prime})}=\big(\tfrac{1}{n}\big)^{\mz(\theta)-c\mz(\theta^\prime)}=
 \big(\tfrac{1}{n}\big)^{\mz(\theta)-\alpha}\to\infty, \quad n\to \infty,
 $$
since $\mz(\theta)<\alpha$. Thus, (\ref{eq:rate_ratio}) holds.

It remains to choose $\kappa$ sufficiently small in order to guarantee
$
2\big(\alpha-\mz(\theta)\big)>c_7\kappa_1^{1+1/\omega_\varsigma}.
$
It yields together with (\ref{eq4:proof-prop1}) that the condition  (\ref{eq:main-assump}) of Theorem \ref{th:generic-adaptive} is fulfilled. Thus, all assumptions of Theorem \ref{th:generic-adaptive} are checked and Proposition \ref{prop:nonexitence-1} follows from this theorem.
\epr

\subsubsection{Proof of Proposition \ref{prop:nonexitence-2}.}
In this proof we will use the second construction presented in Section \ref{sec:sunsec-preliminary_remarks}.

\noindent Let $0<\kappa\leq Q/2$ and $\delta>0$  the constants whose choice will be done later. Set
\begin{eqnarray*}
M&=&\big(\tfrac{\delta\ln(n)}{n^2}\big)^{\frac{\tau_{\varsigma}(q)}{2-2/q-\tau_{\varsigma}(q)}};
\\*[2mm]
\sigma_l&=&\kappa^{\frac{1}{\beta_l}}C_1^{-\frac{1}{\beta_l}}L_l^{-\frac{1}{\beta_l}}
\big[(\kappa/C_1)^{\frac{1}{\beta}}\BL^{-1}M\big]^{\frac{1/\beta_l(1/r_l-1/q)}{\tau_{\varsigma}(q)}};
\\*[2mm]
A&=&\kappa \big[(\kappa/C_1)^{\frac{1}{\beta}}\BL^{-1}M\big]^{-\frac{1/q}{\tau_{\varsigma}(q)}},\quad\; \BL=\prod_{l=1}^d L_l^{\frac{1}{\beta_l}}
\end{eqnarray*}
and remark that
\begin{equation}
 \label{eq5:proof-prop1}
M\bma=\big[(\kappa/C_1)^{\frac{1}{\beta}}\BL^{-1}M\big]^{\frac{1}{\tau_{\varsigma}(q)}}.
\end{equation}
Since $M\to\infty, n\to\infty$ and $\tau_{\varsigma}(q)<0$ one has $M\bma\to 0$ when $n\to\infty$. So, furthermore, we will assume that $M\bma<1$.
Let us  check that $\sigma_l\leq 1$ for any $l=1,\ldots,d$. Indeed, since $M\to \infty, n\to\infty$, $r_l\geq q$ and $\tau_{\varsigma}(q)<0$ we have for all $n$ large enough
 $$
 \sigma_l\leq C_1^{-\frac{1}{\beta_l}}L_l^{-\frac{1}{\beta_l}}\kappa^{\frac{1}{\beta_l}}.
 $$
 Hence, choosing $\kappa$ sufficiently small, guaranteeing
 $
 \kappa\max_{l=1,\ldots,d}C_1^{-1}L_l^{-1}\leq 1
 $
 we come to the required constraints.

 \vskip0.1cm

$\mathbf{1^0}.\;$ First we remark that $x\mapsto A\sum_{m=1}^M y_m\Lambda\big((x-x_{m})/\vec{\sigma}\big)$ belongs to $\bB_q(Q/2)$ since $\kappa\leq Q/2$ and $\|\Lambda\|_q\leq 1$.
Next, we easily check that
the condition
 (\ref{eq:cond-lemma-fromGL14}) of Lemma \ref{lem:lem-from-GL14} holds. It yields together with (\ref{eq:f_0-in-2classes-3})
 $$
 f_y\in\bN_{\vec{r},d}\big(\vec{\beta},\vec{L}\big)\cap\bB_q\big(Q\big),\quad \forall y\in\{-1,1\}^M,\;\forall t\in(0,1).
 $$
 Since (\ref{eq:f_0-in-2classes-3}) also guarantees that $f_{0,t}\in\bN_{\vec{s},d}\big(\vec{\gamma},\vec{L}^\prime\big)\cap\bB_{q^\prime}\big(Q^\prime\big)$ we can assert that Assumption \ref{ass:density_xi_in_class} is verified.

 \vskip0.1cm

$\mathbf{2^0}.\;$ Choose
$
t=M\bma.
$
It yields,
$
\mh_l:=t^{\frac{\rho_l}{\rho}}=\big[(\kappa/C_1)^{\frac{1}{\beta}}\BL^{-1}M\big]^{\frac{\rho_l}{\rho\tau_{\varsigma}(q)}}.
$
Since $\rho_l\leq 1/\beta_l(1/r_l-1/q)$, $\rho\geq 1$, $\kappa\max_{l=1,\ldots,d}C_1^{-1}L_l^{-1}\leq 1$ and $M\to\infty$, $n\to \infty$ and $\tau_{\varsigma}(q)<0$ we obtain that for all $n$ large enough
\begin{equation}
 \label{eq6:proof-prop1}
\sigma_l\leq \mh_l, \quad \forall l=1,\ldots,d.
\end{equation}
Consider rectangle
$
\Pi=:\otimes_{l=1}^d\big[2,2+2\mh_l\big].
$
In view of (\ref{eq6:proof-prop1}) there exist $\mathbf{M}$ disjoint rectangles  with sides of the length  $\sigma_l, l=1,\ldots, d$ belonging to $\Pi$, where
$$
\mathbf{M}:=\prod_{l=1}^d\Big\lfloor\frac{2\mh_{l}}{\sigma_l}\Big\rfloor\geq \frac{\prod_{l=1}^d\mh_{l}}{\prod_{l=1}^d\sigma_l}=M.
$$
Choosing exactly $M$ of these rectangles we come to the collection $\{\Pi_m, m=1,\ldots, M\}$. Thus, Assumption \ref{ass:func_lambda} holds.
Note that
\begin{equation}
\label{eq600:proof-prop1}
f_{0,t}(x)=\mathbf{c}\big[(\kappa/C_1)^{\frac{1}{\beta}}\BL^{-1}M\big]^{-\frac{1/q}{\tau_{\varsigma}(q)}}\mathbf{a}^{d}\mathbf{N}^{-d},\quad \forall x\in \Pi.
\end{equation}
Since $\Pi_m\subset\Pi$ for any $m=1,\ldots,M$, assuming that
$
\kappa<\mathbf{c}\mathbf{a}^{d}\mathbf{N}^{-d}
$
we can assert that Assumptions \ref{ass:density_0} and \ref{ass:one-of-two} are fulfilled. Moreover, for any $m=1,\ldots,M$
\begin{eqnarray}
\label{eq7:proof-prop1}
S_m&=&\int_{\cX} \tfrac{\Lambda^2_m(x)}{f_0(x)}\nu(\rd x)=\mathbf{c}^{-1}\mathbf{a}^{-d}\mathbf{N}^{d}A^{-1}\int_{\bR^d} \Lambda^2_m(x)\rd x
= c_8A\bma \|\Lambda\|_2^2
\nonumber
\\*[2mm]
&=& c_9\kappa\big[(\kappa/C_1)^{\frac{1}{\beta}}\BL^{-1}\big]^{\frac{1-1/q}{\tau_{\varsigma}(q)}}M^{\frac{1-1/q-\tau_{\varsigma}(q)}{\tau_{\varsigma}(q)}}
=c_{10}\kappa^{\frac{\tau_{\varsigma}(1)}{\tau_{\varsigma}(q)}}M^{\frac{1-1/q-\tau_{\varsigma}(q)}{\tau_{\varsigma}(q)}}.
\end{eqnarray}


$\mathbf{3^0}.\;$ Set for brevity $\kappa_1=c_{10}\kappa^{\frac{\tau_{\varsigma}(1)}{\tau_{\varsigma}(q)}}$ and note that $\kappa_1\to\infty$ when $\kappa\to 0$. Note also that
$$
nM^{\frac{1-1/q-\tau_{\varsigma}(q)}{\tau_{\varsigma}(q)}}=n\Big(\tfrac{\delta\ln(n)}{n^2}\Big)^{\frac{1-1/q-\tau_{\varsigma}(q)}{2-2/q-\tau_{\varsigma}(q)}}
=\big(\delta\ln(n)\big)^{\frac{1-1/q-\tau_{\varsigma}(q)}{2-2/q-\tau_{\varsigma}(q)}}n^{\frac{\tau_{\varsigma}(q)}{2-2/q-\tau_{\varsigma}(q)}}\to 0,\;\;n\to\infty,
$$
since $\tau_{\varsigma}(q)<0$.
Hence  $nS_m\to 0, n\to\infty$, in view of (\ref{eq7:proof-prop1})  we get
\begin{eqnarray}
\label{eq8:proof-prop1}
\prod_{m=1}^M \int_{\bR^2}e^{\mathbf{i} nS_m uv}\bP(\rd u)\bP(\rd v)&\leq& \Big( 1+c_{11}\kappa_1^2n^2M^{\frac{2-2/q-2\tau_{\varsigma}(q)}{\tau_{\varsigma}(q)}}\Big)^M
\leq e^{c_{11}\kappa_1^2n^2M^{\frac{2-2/q-\tau_{\varsigma}(q)}{\tau_{\varsigma}(q)}}}
\nonumber
\\*[2mm]
&=&e^{c_{11}\kappa_1^2\delta\ln(n)}=n^{c_{11}\kappa_1^2\delta}.
\end{eqnarray}
Let $\theta^{\prime\prime}\in \Theta^{\prime\prime}[\theta]$ and  $\mz(\theta)<\alpha< \mz(\theta^{\prime\prime})$ be fixed. Set $c=\frac{\alpha}{\mz(\theta^{\prime\prime})}$ and note that $c<1$. On the other hand
 $$
 \alpha_n(c):=\tfrac{\phi_n(\cF_\vartheta)}{\phi^{c}_n(\cF_{\vartheta^{\prime\prime}})}=\big(\tfrac{1}{n}\big)^{\mz(\theta)-c\mz(\theta^{\prime\prime})}=
 \big(\tfrac{1}{n}\big)^{\mz(\theta)-\alpha}\to\infty, \quad n\to \infty,
 $$
since $\mz(\theta)<\alpha$. Thus, (\ref{eq:rate_ratio}) holds. It remains to choose $\delta$ sufficiently small in order to guarantee
$$
2\big(\alpha-\mz(\theta)\big)>c_{11}\kappa_1^2\delta.
$$
It yields together with (\ref{eq8:proof-prop1}) that the condition  (\ref{eq:main-assump}) of Theorem \ref{th:generic-adaptive} is fulfilled.

\vskip0.1cm

$\mathbf{4^0}.\;$ Note that for all $n$ large enough
$$
\big|\|f_y\|_2-\|f_{0,t}\|_2\big|\geq c_{12}\big|\|G_t+F_y\|_2-\|G_{t}\|_2\big|,
$$
where we have put $F_y(x)=A\sum_{m=1}^M y_m\Lambda\big((x-x_{m})/\vec{\sigma}\big)$. Recall also that $t=M\bma$.

Taking into account (\ref{eq600:proof-prop1}) and that $\int \Lambda =0$ we obtain
\begin{eqnarray*}
\big|\|f_y\|_2-\|f_{0,t}\|_2\big|&\geq&\frac{ c_{13}A^2M\bma}{2\|G_{t}\|_2+A\sqrt{M\bma}}\geq\frac{ c_{14}A^2M\bma}{(M\bma)^{1/2-1/q}+A\sqrt{M\bma}}
=\frac{ c_{14}A^2\sqrt{M\bma}}{(M\bma)^{-1/q}+A}
\\*[2mm]
&=& c_{15}M^{\frac{1/2-1/q}{\tau_{\varsigma}(q)}}=c_{16}\Big(\tfrac{\sqrt{\ln(n)}}{n}\Big)^{\frac{1-2/q}{2-2/q-\tau_{\varsigma}(q)}}
\end{eqnarray*}
We assert that Assumption \ref{ass:on-Phi} is fulfilled with
$
\psi_n=2^{-1}c_{16}\Big(\tfrac{\sqrt{\ln(n)}}{n}\Big)^{\frac{1-2/q}{2-2/q-\tau_{\varsigma}(q)}}.
$
Thus, all assumptions of Theorem \ref{th:generic-adaptive} are checked and Proposition \ref{prop:nonexitence-2} follows from this theorem.
\epr

\section{Appendix.}

\paragraph{Proof of Lemma \ref{lem:bound_for_f_y}.}

As we already mentioned $f^*_{y}\equiv f_{y}$ if $\lambda_m=0$ for all $m=1,\ldots, M$ and, therefore the assertion of the lemma is true,
since in view of Assumption \ref{ass:one-of-two} $f_{y}$ is probability density for all $y\in\bY^M$ in this case.

\noindent In the remaining case $\bY\subseteq [0,1]^M$  and $\Lambda:\cX\to\bR_+$ for all $m=1,\ldots, M$.
We obviously have
$$
(1-\rho_y)=(1-\Sigma_M)\Big(1-\tfrac{\rho_y-\Sigma_M}{1-\Sigma_M}\Big).
$$
Note that in view of (\ref{eq:ass-on-parameters}) for any $y\in\bY_\rho$
$$
\Big|\tfrac{\rho_y-\Sigma_M}{1-\Sigma_M}\Big|\leq 16\cD_M <n^{-\frac{1}{2}}.
$$
Applying the  inequality $-t\geq \ln(1-t)\geq -t-t^2$ which holds for all $t>0$ small enough,  we assert
$$
e^{-\frac{\rho_y-\Sigma_M}{1-\Sigma_M}}\geq \Big(1-\tfrac{\rho_y-\Sigma_M}{1-\Sigma_M}\Big)\geq e^{-\frac{\rho_y-\Sigma_M}{1-\Sigma_M}-\left(\frac{\rho_y-\Sigma_M}{1-\Sigma_M}\right)^2}
\geq  e^{-\frac{\rho_y-\Sigma_M}{1-\Sigma_M}-\frac{1}{n}}
$$
for all $n$ large enough.
Since  $1-\Sigma_M\geq 0.75$ in view of Assumption \ref{ass:one-of-two}, we have for any $y\in\bY_\rho$
\begin{align*}
(1-\Sigma_M)e^{-\frac{\rho_y-\Sigma_M}{1-\Sigma_M}}f_0(x)&+\sum_{m=1}^My_m\Lambda_m(x)
\geq f_{y}(x)
\\
&\geq e^{-\frac{1}{n}}(1-\Sigma_M)e^{-\frac{\rho_y-\Sigma_M}{1-\Sigma_M}}f_0(x)+\sum_{m=1}^My_m\Lambda_m(x)
\\
&\geq e^{-\frac{1}{n}}\bigg[(1-\Sigma_M)e^{-\frac{\rho_y-\Sigma_M}{1-\Sigma_M}}f_0(x)+\sum_{m=1}^My_m\Lambda_m(x)\bigg].
\end{align*}
Here we have used that $f_0\geq 0$ because it is a probability density. To get the last inequality we have also used that
$y\in[0,1]^M$  and $\Lambda:\cX\to\bR_+$ for all $m=1,\ldots, M$.
The lemma is proved.
\epr

\paragraph{Proof of Lemma \ref{lem:probab}.}

Since the random variable  $|\chi_K|$ takes  values in $[0,2TKA_K]$ we have
\begin{align*}
W_{J,K}(b) = 1 + 2Jb \int_0^{2TKA_K} y(1+by^2)^{J-1} P\big(
 |\chi_K|\geq y\big) \rd y.
\end{align*}
Set
$V^2_K=E(\zeta_1^2)\sum_{k=1}^K a^2_k$
 and  note that
that $\chi$ is a sum of i.i.d. centered random variables taking values in $[-2TA_K,2TA_K]$.
Therefore applying the Bernstein  inequality together with trivial inequality $1+z\leq e^z$ we obtain
 \begin{align}
 \label{eq:ol00}
 W_{J,K}(b) &\leq 1+2Jb\int_0^{2TKA_K} y \exp\Big\{ b(J-1)y^2 - \tfrac{y^2}{2V^2_K +\frac{4}{3}TA_Ky}\Big\}\rd y
 \nonumber
 \\*[2mm]
 &\leq 1+2 Jb\int_0^{2TKA_K} y \exp\Big\{bJy^2 - \tfrac{y^2}{2V^2_K +\frac{8}{3}T^2A^2_KK}\Big\}\rd y
 \nonumber
 \\*[2mm]
 &\leq1+2 Jb\int_0^{\infty} y \exp\Big\{bJy^2 - \tfrac{y^2}{D_K}\Big\}\rd y.
 \end{align}
Here we have also used that $V^2_K\leq T^2A^2_KK$.
 Setting
$
 \omega=D_K^{-1}-bJ
$
and noting that $\omega>0$ we get
\begin{align*}
\label{eq:ol001}
 W_{J,K}(b) \leq 1+2Jb\int_0^{\infty} y e^{-\omega y^2}\rd y=1+Jb\omega^{-1}\leq 2.
 \end{align*}
in view of the condition of the lemma. This completes the proof.
\epr


\bibliographystyle{agsm}

\end{document}